\documentclass[12pt,a4paper]{amsart}
\usepackage{amsfonts}
\usepackage{amsthm}
\usepackage{amsmath}
\usepackage{amscd}
\usepackage[latin2]{inputenc}
\usepackage{t1enc}
\usepackage[mathscr]{eucal}
\usepackage{indentfirst}
\usepackage{graphicx}
\usepackage{graphics}
\usepackage{pict2e}
\usepackage{epic}
\numberwithin{equation}{section}
\usepackage[margin=2.9cm]{geometry}
\usepackage{epstopdf} 
\usepackage{esvect}
\usepackage{mathtools}
\usepackage[english]{babel}
\usepackage{xspace}
 
\usepackage[english]{babel}
\usepackage{verbatim}

 \newtheorem{theorem}{Theorem}[section]
 
 \newtheorem{lemma}{Lemma}[section]

 \numberwithin{equation}{section}

\theoremstyle{definition}
\newtheorem{definition}{Definition}[section]

\usepackage{enumitem}

\begin{document}

\setlength{\abovedisplayskip}{3.78pt}%
\setlength{\belowdisplayskip}{3.78pt}%
\setlength{\abovedisplayshortskip}{3.78pt}%
\setlength{\belowdisplayshortskip}{3.78pt}%
\setlist[itemize]{noitemsep, topsep=0pt}
\allowdisplaybreaks

\title[Global Solution to a Non-linear Wave Equation of Liquid Crystal]{
Global Solution to a Non-linear Wave Equation of Liquid Crystal in the Constant Electric Field
}

\author{ Linjun Huang }

\begin{abstract} We construct a global conservative weak solution to the Cauchy problem for the non-linear variational wave equation $v_{tt} - c(v)(c(v)v_x)_x + \frac{1}{2}(v+v^3)= 0$ where $c(\cdot)$ is any smooth function with uniformly positive bounded value. This wave equation is derived from a wave system modelling nematic liquid crystals in a constant electric field.
\end{abstract}

\maketitle

\section{ \textbf{Introduction}} 
\indent In this paper, we study a wave equation modelling the nematic liquid crystal in one space dimension with electric field applied. In the nematic phase, the orientation of the molecules can be described by a field of unit vector $\mathbf{n} (x,t) \in S^2$, the unit sphere. The famous Oseen-Frank potential energy density $W$ associated with the director field $\mathbf{n}$ is defined by
 $$ W(\mathbf{n}, \nabla \mathbf{n}  ) = \alpha |\mathbf{n} \times (\nabla \times \mathbf{n})|^2 + \beta (\nabla \cdot \mathbf{n})^2 + \gamma (\mathbf{n} \cdot \nabla \times \mathbf{n})^2,$$
where $\alpha, \beta$ and $\gamma$ are positive elastic constants of the liquid crystal. $\alpha$ represents the splay phenomenon of the nematic liquid crystal, $\beta$ represents the bend phenomenon, and $\gamma$ represents the twist phenomenon. When the kinetic energy are neglected in studies of nematic liquid crystals, by variational principle, we obtain an elliptic partial differential equation \cite{newnewname1}. When we include the kinetic energy on modelling the nematic liquid crystal in one space dimension without any fields applied, we can formulate it as a non-linear wave equation 
 \begin{equation} \label{ex:0101}
  \displaystyle u_{tt} - c(u)[c(u) u_x]_x = 0 ,
 \end{equation}
with smooth function $u$.  \\
 \indent We study the nematic liquid crystal under the a constant electric field with the electric energy density described by
$$ \displaystyle f_{electric} = - \frac{1}{2} \mathbf{P} \cdot \mathbf{E} =\frac{1}{2}\varphi \mathbf{E}^2 + \frac{1}{2} \eta( \mathbf{E} \cdot \mathbf{n})^2 ,$$ 
 where $\mathbf{P}$ is the polarization, $\mathbf{E}$ is the electric field. We assume that the applied field is neither parallel nor perpendicular to $\mathbf{n}$. $\varphi$ and $\eta$ are positive constants related to permittivity and dielectric constants \cite{name6}.

\subsection{Known results} $ $\\
\indent For the equation (\ref{ex:0101}), Glassey, Hunter, and Zheng \cite{name18} showed that the smooth solutions develop singularities in finite time. Also, Zhang and Zheng \cite{name15} studied that under weak conditions on the initial data which allow the solutions to have blow-up singularities and they established approximate solutions with estimates along precompactness using Young measure methods.\\
\indent Our main reference is \cite{name5}. For the Cauchy problem for (\ref{ex:0101}) with initial data $u(0,x) = u_0 (x)$, $u_t(0,x) = u_1(x)$, Bressan and Zheng \cite{name5} proved the existence of a conservative weak solution by method of characteristics. They constructed conservative weak solution by introducing new sets of dependent and independent variables and showed that the solution can be obtained as the fixed point of a contraction transformation. See also \cite{newname1}. Compared with \cite{name5}, our energy equation has new terms from the applied electric field. These terms can be expressed as $v^2+v^4$. To solve this problem, we need to do some modification on the proof in \cite{name5} based on the observation that $v \in H^1 \cap L^4$ and $v^2$, $v^4$ are the lower order terms in the energy equation.\\
\indent For the Cauchy problem for (\ref{ex:0101}) with initial data, Bressan, Chen and Zhang \cite{name3} proved the uniqueness of conservative solutions. Brassan and Huang \cite{name4} constructed dissipative solutions relying on Kolmogorov's compactness theorem. Zhang and Zheng studied the existence and regularity properties of classical and weak solutions using the Young measure theory in \cite{name14} and proved the global existence of weak solutions in \cite{name16}. For $C^3$ initial data, Bressan and Chen \cite{name1} showed that the conservative solutions are piecewise smooth in $t$-$x$ plane. In \cite{name2}, Bressan and Chen constructed a metric that renders the flow uniformly Lipschitz continuous on bounded subsets of $H^1(\mathbb{R})$. Zhang and Zheng \cite{name17} studied the existence of global weak solutions to the initial value problem (\ref{ex:0101}) with general initial data $(u(0),u_t(0))=(u_0,u_1) \in W^{1,2} \times L^2$ with wave speed satisfying $c'(\cdot) \geq 0$ and $c'(u_0(\cdot)) > 0$.\\
\indent For a wave system modelling nematic liquid crystals in one space dimension, Chen and Zheng \cite{name7} studied the global existence and singularity formation. Huang and Zheng \cite{name8} established the global existence of smooth solutions. Zhang and Zheng \cite{name9} constructed a weak global solutions to the Cauchy problem for a system of two variational wave equations on the real line and \cite{name10} showed the global weak solutions to the initial value problem for a complete system of variational wave equations modelling liquid crystals in one space dimension. \\
\indent Equation (\ref{ex:0101}) has asymptotic uni-directional wave equation
\begin{equation} \label{neweq:115}
(u_t + u^n u_x)_x = \frac{1}{2} n u^{n-1} (u_x)^2,
\end{equation}
derived by Hunter and Saxton via weakly nonlinear geometric optics. In \cite{name11}-\cite{name13}, Zhang and Zheng studied the global existence, uniqueness, and regularity of the dissipative and conservative solutions to (\ref{neweq:115})$(n=1,2)$ with $L^2$ initial data. 
\subsection{Main theorems} $ $ \\
\indent  Our main results are stated as follows.
For the nematic liquid crystal under electric field, we obtain a Cauchy problem
\begin{equation} \label{eq:011}
\displaystyle v_{tt} - c(v)(c(v)v_x)_x + \frac{1}{2}(v+v^3)= 0 ,
\end{equation}
with the initial data
\begin{equation} \label{eq:012}
  v(0,x)=v_0(x), \,\,\,\, v_t(0,x)=v_1(x).
\end{equation}
\indent For the smooth function $c(\cdot)$, we assume that $c: \mathbb{R}  \mapsto \mathbb{R^{+}} $ is a bounded and uniformly positive function.
\begin{definition}
The definition of weak solution.\\
\indent We say that for all test function  $ \phi \in C^1_c$, the function $ \displaystyle v \in H^1 \cap L^4$ satisfies the following integral: 
\begin{equation} \label{eq:014}
\iint \phi_t v_t - [c(v) \phi ]_x [c(v) v_x] - \frac{\phi}{2} (v + v^3) dx dt = 0,
\end{equation}
is a weak solution to the equation (\ref{eq:011}).
\end{definition}
\begin{definition}
The definition of energy conservative weak solution.\\
\indent For $v_1$ and $v_0$ defined in (\ref{eq:012}), we define the ground state energy $ \mathcal{E}_0 $ as:
\begin{equation} \label{eq:016}
\mathcal{E}_0 \coloneqq \frac{1}{2} \int \Big\{  v^2_1(x) + c^2 (v_0(x)) [v_0^2(x)]_x + \frac{v_0^2(x)}{2} + \frac{v^4_0(x)}{4} \Big\} dx.
\end{equation}
The function $ \displaystyle v \in H^1 \cap L^4$ is a energy conservative weak solution if it satisfying
\begin{equation} \label{eq:015}
  \mathcal{E} (t) \coloneqq \frac{1}{2} \int \Big\{  v^2_t(t,x) + c^2 (v(t,x)) v_x^2(t,x) + \frac{v^2(t,x)}{2} + \frac{v^4(t,x)}{4} \Big\} dx = \mathcal{E}_0,
\end{equation}
for almost every $t \in \mathbb{R}$.
\end{definition}
\begin{theorem} \label{thm_main}
Assume that c: $\mathbb{R} \mapsto [\mathcal{K}^{-1}, \mathcal{K}]$ is a smooth function for some $ \mathcal{K} > 1 $. $ v_0(x)$ and $v_1(x)$ are stated in (\ref{eq:012}). Also assume that the initial data $v_0(x) $ is absolutely continuous, $(v_0(x))_x \in H^1 \cap L^4$, and $v_1(x) \in H^1 \cap L^4 $. Then (\ref{eq:011})-(\ref{eq:012}) can be considered as a Cauchy problem admitting a weak solution $v(t,x)$ defined for all $(t,x) \in \mathbb{R} \times \mathbb{R}$. Moreover, in the $t$-$x$ plane, $v(x,t)$ is locally H\"older- $\frac{1}{2}$ continuous. For all $ 1 \leq p < 2$, the map $ t \mapsto v(t, \cdot) $ is continuously differentiable with values in $L^p_{loc}$. The weak solution $v(t. \cdot)$ is Lipschitz continuous with respect to $L^2$ distance. So, for all $t,s \in \mathbb{R}$,
\begin{equation} \label{eq:013}
   \| v(t, \cdot) - v(s, \cdot) \|_{L^2} \leq L |t - s|.
\end{equation}
For all test function  $ \phi \in C^1_c$, the equation (\ref{eq:011}) satisfies (\ref{eq:014}).
\end{theorem}
\begin{theorem} \label{sec}
 A family of weak solutions to the Cauchy problem (\ref{eq:011})-(\ref{eq:012}) can be obtained with the properties:
 \begin{equation} \label{eq:017}
  \mathcal{E} (t) \leq \mathcal{E}_0 .
 \end{equation}
Let a sequence of initial condition satisfies:
$$ \| (v^n_0(x))_x - (v_0(x))_x \|_{L^2} \to 0 ,$$
$$ \| v^n_1(x) - v_1(x) \|_{L^2} \to 0 .$$ 
 Also, $ u^n \to u$ uniformly on bounded subsets of the $t$-$x$ plane and $v^n_0 \to v_0$ on  compact sets as $n \to \infty$.
\end{theorem}
\begin{theorem} \label{thr}
 There exists a continuous family of positive Radon measures $ \{ \mu_t : t \in \mathbb{R} \} $ . This family of positive Radon measure is defined on the real line and it satisfies the following properties: \\
 (i) $ \mu_t (\mathbb{R}) = \mathcal{E}_0$ for any time $t$. \\
 (ii) With respect to Lebesgue measure, the absolutely continuous part of $\mu_t$ has density $\displaystyle \frac{1}{2}( v_t^2 + c^2(v) v_x^2 + \frac{1}{2} v^2 + \frac{1}{4} v^4 )$. \\
 (iii) The singular part of $\mu_t$ has measure zero on the set where $ c'(v) = 0$. 
\end{theorem}
The paper is organized as follows. In section 2, we derive the energy equation and introduce a new set of dependent variables. Based on those dependent variables, we formulate a set of equations in terms of the new variables. This set of equations is equivalent to (\ref{eq:011}). In section 3, we use a transformation in a Banach space. In the transformation, we find the suitable weighted norm. This shows that there is a unique solution to the set of equations in terms of the new variables. In section 4, we show that the integral (\ref{eq:014}) holds  and the H\"older-$\frac{1}{2}$ continuous condition holds. In section 5, we show that (\ref{eq:017}) holds and the Lipschitz condition on the map $ t \mapsto v(t, \cdot)$ and provide a proof of Theorem 1.2. On section 6, we study the maps of $t \mapsto u_x (t, \cdot)$ and $ t \mapsto u_t(t, \cdot)$, and complete the proof of Theorem \ref{thm_main}. We provide a proof of Theorem 1.3 in section 7.
\section{\textbf{Variable Transformations}} 
\subsection{Derivation of (\ref{eq:011}) } 
Equation (\ref{eq:011}) has some physical origins. In the context of nematic liquid crystals, we introduce the famous Oseen-Frank potential energy density $W$ is given by
\begin{equation} \label{eq:18}
 W( \mathbf{n},\nabla \mathbf{n} ) = \alpha |\mathbf{n} \times (\nabla \times \mathbf{n})|^2 + \beta(\nabla \cdot \mathbf{n})^2 + \gamma(\mathbf{n} \cdot \nabla \times \mathbf{n})^2 .
 \end{equation}
 As stated in \cite{name6}, in a electric field, the electric energy of the liquid crystal per unit volume is given by
\begin{equation} \label{eq:19}
 f_{electric} = - \frac{1}{2} \mathbf{P} \cdot \mathbf{E} = \frac{1}{2} \varphi \mathbf{E}^2 + \frac{1}{2} \eta ( \mathbf{E} \cdot \mathbf{n})^2 .
 \end{equation}

\indent We discuss that when the electric energy is low when the applied electric field is normal to the liquid crystal director. And $\varphi$ and $\eta$ are some positive constants related to the permittivity so that (\ref{eq:19}) is equivalent to $-\frac{1}{2} | \mathbf{n} \cdot \mathbf{E}^{\perp} |^2 + 1$. And we denote $\mathbf{E}^{\perp}$ as a vector such that $\mathbf{E} \cdot \mathbf{E}^{\perp} = 0$. So, the electric energy can be described as $-\frac{1}{2} | \mathbf{n} \cdot \mathbf{E}^{\perp} |^2$. By the property of the potential energy, the action can be describe as $ \frac{1}{2} |\mathbf{n}_t|^2 - W(\mathbf{n}, \nabla \mathbf{n} ) + \frac{1}{2} |\mathbf{n} \cdot \mathbf{E}^{\perp}|^2$. We denote that $\mathbf{F} = \nabla \mathbf{n}$, so that $W(\mathbf{n}, \nabla \mathbf{n}) = W(\mathbf{n}, \mathbf{F})$. And we define `$:$' as for $\displaystyle{ A,B \in \mathbb{R}^{n \times n}}$, $\displaystyle{ A:B = \sum\limits^n_{j=1} \sum\limits^n_{i=1} A_{ij}B_{ji}}$. The below Euler-Largrangian equation derived from the least action principle 
\begin{align*}
 &L(\mathbf{n}) = |\mathbf{n}_t|^2 - 2 W(\mathbf{n}, \mathbf{F}) + |\mathbf{n} \cdot \mathbf{E}^{\perp}|^2 \\
& \tilde{L} (\mathbf{n}, \lambda) = L(\mathbf{n}) - \lambda \cdot (|\mathbf{n}|^2 - 1)\\
& \mathcal{L} (\mathbf{n}, \lambda) = \iint |\mathbf{n}_t|^2 - 2 W(\mathbf{n}, \mathbf{F} ) + |\mathbf{n} \cdot \mathbf{E}^{\perp}|^2  - \lambda (|\mathbf{n}|^2 - 1) dx dt 
\end{align*}
From the least action principle,
\begin{align*}
 0 = \delta \mathcal{L} (\mathbf{n}, \lambda) = &\iint \mathbf{n}_t \delta n_t - (\partial_n W \cdot \delta n + \partial_F W : \delta \nabla \mathbf{n}) + (\mathbf{n} \cdot \mathbf{E}^{\perp}) \delta n \cdot \mathbf{E}^{\perp} - \lambda \mathbf{n} \delta n \\
& -\frac{1}{2} \delta \lambda (|\mathbf{n}|^2 - 1) \,\, dx dt
\end{align*}
The above equation can be describe as follow: 
\begin{equation} \label{eq:110}
 \mathbf{n}_{tt} -  \partial_n W( \mathbf{n}, \partial_x \mathbf{n} ) - \partial_x [ \partial_{\partial_x} W(\mathbf{n}, \partial_x \mathbf{n})]+  (\mathbf{n} \cdot \mathbf{E}^{\perp}) \cdot \mathbf{E}^{\perp} =  \lambda \mathbf{n}, \,\,\,\, \mathbf{n} \cdot \mathbf{n} = 1.
\end{equation}

For planar deformations depending on a single space variable $x$, the director field has the from $\mathbf{n} = \cos  u(x,t) \mathbf{i} + \sin u(x,t) \mathbf{j}$. The dependent variable $u \in \mathbb{R}^1$ measures the angle of the director field to the x direction. In this case, we have the wave speed $c$ given specifically by $ c^2(u) = \gamma \cos  ^2 u + \alpha \sin ^2 u$ and $ \alpha = \beta$, $\gamma = 0$ for $\delta W(\mathbf{n}, \nabla \mathbf{n} )$.
In one space dimension, (\ref{eq:110}) becomes,
\begin{equation} \label{eq:111}
\lambda n_i  = \partial_{tt} n_i + \partial_{n_i} W( \mathbf{n}, \partial_x \mathbf{n}) - \partial_x[\partial_{\partial_x n_i} W(\mathbf{n}, \partial_x \mathbf{n})] + (\mathbf{n} \cdot \mathbf{E}^{\perp}) \cdot \mathbf{E}^{\perp} , \,\,\, 
 \end{equation}
 for  i=1,2,3.
\begin{equation} \label{eq:112}
\lambda = - | \mathbf{n}_t|^2 +(\beta + 2(\gamma - \beta)n_1^2)|\partial_x \mathbf{n}|^2 +(\gamma - \beta ) n_1  \partial_x^2 n_1 - |\mathbf{n} \cdot \mathbf{E}^{\perp}|^2.
\end{equation}
In 1-d case, the Oseen-Frank potential energy states in(\ref{eq:18}) is
\begin{equation} \label{eq:113}
 W(\mathbf{n}, \partial_x \mathbf{n}) = \frac{\alpha}{2} (\partial_x n_1)^2 + \frac{\beta}{2} [(\partial_x n_2)^2 + (\partial_x n_3)^2] + \frac{1}{2} (\gamma - \beta ) n_1^2 |\partial_x \mathbf{n}|^2.
\end{equation}
We denote that $ c^2_1 = \alpha + (\gamma - \alpha) n_1^2$ and $ c^2_2 = \beta + (\gamma - \beta) n_1^2$ . With (\ref{eq:112}) and (\ref{eq:113}), we can compute that the left hand side of (\ref{eq:111}) becomes
\begin{equation} \label{eq:114}
 \partial_{tt} n_1 - \partial_x[c_1^2(n_1) \partial_x n_1] + (\mathbf{n} \cdot \mathbf{E}^{\perp}) \delta \mathbf{n} \cdot E^{\perp} = \{ -|\mathbf{n}_t|^2 + (2 c_2^2 - \gamma)|\mathbf{n}_x|^2 + 2( \alpha - \beta)(\partial_x n_1)^2 - |\mathbf{n} \cdot \mathbf{E}^{\perp}|^2 \} n_1.
\end{equation}
In particular, taking $\alpha = \beta$ in (\ref{eq:113}), we let 
$$ c^2(u) = c_1^2(n_1) = c_2^2(n_1).$$
By plug in  $\mathbf{n} = \cos  u(x,t) \mathbf{i} + \sin u(x,t) \mathbf{j}$ into (\ref{eq:114}).
Thus, we plug in the above calculations into (\ref{eq:114}) and get
\begin{equation} \label{eq:115}
\begin{split}
- u_{tt} + c^2(u) u_{xx} + (- \gamma \cos (u) \sin (u) + \alpha \cos (u) \sin (u)) u^2_x  + \frac{\cos ^2(u)}{\sin (u)} + \frac{\cos (u)}{\sin (u)} = 0 ,  \\
  u_{tt} - c(u) ( c(u) u_x )_x  - \frac{\cos (u) (1 + \cos (u) )}{\sin (u)} = 0.
\end{split}
\end{equation}
For the asymptotic equation (\ref{eq:115}), we do the Taylor expansion for $\displaystyle{ \frac{\cos (u) (1 + \cos (u) )}{\sin (u)} }$ around point $ u = \pi$ ignoring some high order terms. So, (\ref{eq:115}) becomes 
$$   u_{tt} - c(u) ( c(u) u_x )_x  + \frac{(u + \pi)^3}{2} +\frac{(u + \pi)}{2} = 0.  $$
We let $ v = u + \pi$ that give us  
\begin{equation} \label{eq:116}
  v_{tt} - c(v) ( c(v) v_x )_x  + \frac{1} {2} (v + v^3) = 0.
\end{equation}
\subsection{Derivation of the energy equation} $ $ \\
\indent From (\ref{eq:116}), we can compute that
 \begin{equation} \label{eq:021}
 \begin{split}
 \int v_t v_{tt} - v_t c(v) [c(v) v_x]_x + \frac{v}{2} v_t + \frac{v^3}{2} v_t dx  &= 0\\
\int (\frac{1}{2} v_t^2)_t + (\frac{c^2(v) v_x^2}{2})_t + (\frac{1}{4} v^2)_t + (\frac{1}{8} v^4)_t dx &= 0.
\end{split}
 \end{equation}
 And from (\ref{eq:021}), the energy equation can be described as
 \begin{equation} \label{eq:022}
  E \coloneqq \frac{1}{2} \left( v_t^2 + c^2(v) v_x^2 + \frac{1}{2} v^2 + \frac{1}{4} v^4 \right).
 \end{equation}
\subsection{Variables transform}
In this section we derive identities that holds for smooth solutions. We first denote variables:
\begin{equation} \label{eq:023}
\left\{
       \begin{array}{ll}
        R \coloneqq v_t + c(v) v_x ,\\
        S \coloneqq v_t - c(v) v_x .
        \end{array}  
        \right.
\end{equation}
Thus, we can write $v_t$ and $v_x$ as follows
\begin{equation} \label{eq:024}
\left\{
       \begin{array}{ll}
        \displaystyle v_t = \frac{R+S}{2} , \\[10pt]
        \displaystyle v_x = \frac{R-S}{2c}.
        \end{array}  
        \right.
\end{equation}
By (\ref{eq:011}), the following identities are valid :
\begin{equation} \label{eq:025}
\left\{
       \begin{array}{ll}
       \displaystyle  S_t +c S_x = \frac{c'}{4c} (S^2 - R^2) - \frac{1}{2} (v + v^3), \\[10pt]
       \displaystyle  R_t -c R_x = \frac{c'}{4c} (R^2 - S^2) - \frac{1}{2} (v + v^3),
        \end{array}  
        \right.
\end{equation}
by the following calculation
\begin{align*}
R_t - c R_x &= (v_t + c v_x)_t - c ( v_t + c v_x)_x \\
            &= - \frac{v}{2} - \frac{v^3}{2} + \frac{c'}{4c} (R^2 - S^2).
\end{align*}
We can compute $ S_t +c S_x $ in the similar way to get (\ref{eq:025}) and denote energy and momentum as 
\begin{equation} \label{eq:026}
 E \coloneqq  \frac{1}{2} \left( v_t^2 + c^2(v) v_x^2 + \frac{1}{2} v^2 + \frac{1}{4} v^4 \right) = \frac{R^2 + S^2}{4} + \frac{1}{8} (2 v^2 + v^4),
\end{equation}
\begin{equation} \label{eq:027}
M \coloneqq - v_t v_x = \frac{S^2 - R^2}{4c}.
\end{equation}
\indent The analysis of (\ref{eq:011}) has a main difficult that the possible breakdown of the regularity solutions. The quantities $v_x$ and $v_t$ can blow up in finite time even with smooth initial data. Thus we need to introduce a new set of dependent variables to deal with the possible unbounded value $R$ and $S$: 
\begin{equation} \label{eq:028}
w \coloneqq 2 \arctan R , \, \,\,\,\,\,\,\,\,\,\, z \coloneqq 2 \arctan  S.
\end{equation}
Thus 
\begin{equation}\label{eq:029}
R = \tan \left(\frac{w}{2}\right), \,\,\,\,\,\,\,\,\,\,\,\,\, S = \tan \left(\frac{z}{2}\right).
\end{equation}
By (\ref{eq:025}), 
\begin{equation} \label{eq:0210}
 \displaystyle w_t - c w_x = \frac{2}{1+R^2} (R_t - c R_x) = \frac{c'}{2c} \frac{R^2 - S^2}{1 + R^2} - \frac{1}{1 + R^2} (v + v^3),
\end{equation}
\begin{equation} \label{eq:0211}
 \displaystyle z_t + c z_x = \frac{2}{1+S^2} (S_t + c S_x) = \frac{c'}{2c} \frac{S^2 - R^2}{1 + S^2} - \frac{1}{1 + S^2} (v + v^3).
\end{equation}
\indent In order to reduce the equation to a semi-linear system, we need to have a further change of variables. The forward characteristics equation and the backward characteristics equation:
$$ \dot{x}^+ = c(v), \,\,\,\,\,\,\,\, \dot{x}^- = c(v). $$
And we denote the characteristics lines pass through the point $(t,x)$ as 
$$ s \to x^+ (s,t,x), \,\,\,\,\,\,\,\,  s \to x^- (s,t,x)  .$$
\begin{figure}
  \centering
      \includegraphics[scale=0.7]{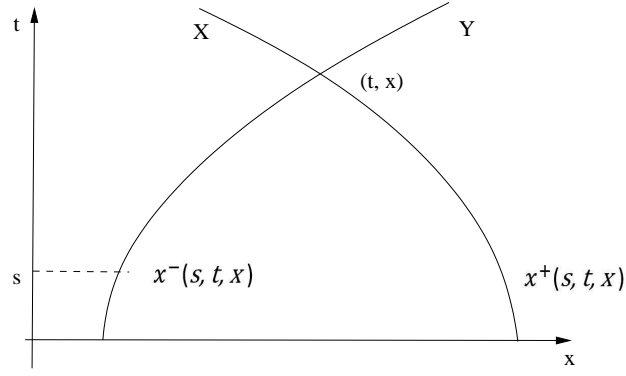}
      \caption{The characteristic curves.}
\end{figure}
So we can use a new coordinate system $(X,Y)$ to represent point $(t,x)$ by
\begin{equation} \label{eq:0212}
 X \coloneqq \int^{x^- (0,t,x)}_0 (1 + R^2(0,x))dx,
 \end{equation}
 \begin{equation} \label{eq:0213}
 Y \coloneqq \int_{x^+ (0,t,x)}^0 (1 + S^2(0,x))dx.
\end{equation} 
(\ref{eq:0212}) and (\ref{eq:0213}) implies that 
\begin{equation} \label{eq:0214}
X_t - c(v) X_x = 0, \,\,\,\,\,\,\,\,\,\,\,\,\,\, Y_t + c(v) Y_x = 0,
\end{equation}
\begin{equation} \label{eq:0215}
(X_x)_t - (c(v) X_x)_x = 0, \,\,\,\,\,\,\,\,\,\,\,\,\,\, (Y_x)_t + (c(v) Y_x)_x = 0.
\end{equation}
Thus, given any smooth function $f$, by using (\ref{eq:0214}),
\begin{equation} \label{eq:0216}
\begin{array}{ll}
        f_t + c(v) f_x = 2c(v) X_x f_X, \\
        f_t - c(v) f_x = 2c(v) Y_x f_Y.
        \end{array}  
\end{equation}
From (\ref{eq:0214}), $ X_t + c(v) X_x = 2c(v)  X_x$. To get (\ref{eq:0216}) we compute directly
\begin{align*}
        f_t + c(v) f_x = f_X X_t + f_Y Y_t + c(v) f_X X_x + c(v)f_Y Y_x = (X_t + c(v) X_x) f_X = 2c(v) X_x f_X, \\
        f_t - c(v) f_x = f_X X_t + f_Y Y_t - c(v) f_X X_x - c(v) f_Y Y_x = (Y_t - c(v) Y_x) f_Y = 2c(v) Y_x f_Y.
\end{align*}
Introducing new variables
\begin{equation} \label{eq:0217}
 p \coloneqq \frac{1 + R^2}{X_x}, \,\,\,\,\,\,\,\,\,\,\,\, q \coloneqq \frac{1 + S^2}{- Y_x}.
\end{equation}
From (\ref{eq:0217}),
\begin{equation} \label{eq:0218}
       \begin{array}{ll}
       \displaystyle{ \frac{1}{X_x} = \frac{p}{1 +R^2} = p \cos ^2(\frac{w}{2}) = \frac{p(1+\cos  w )}{2},}\\
       \displaystyle{ \frac{1}{-Y_x} = \frac{q}{1 +S^2} = q\cos ^2(\frac{z}{2}) = \frac{q(1+\cos  z )}{2}.}
        \end{array}  
\end{equation}
By applying (\ref{eq:0210})-(\ref{eq:0211}) to (\ref{eq:0211}), 
$$ \displaystyle w_t - c w_x = 2 c \frac{1+S^2}{q} w_Y = \frac{c'}{2c} \frac{R^2 - S^2}{1+R^2} + \frac{1}{1+R^2} (-v-v^3),$$
$$ \displaystyle z_t + c z_x = 2 c \frac{1+ R^2}{p} z_x = \frac{c'}{2c} \frac{S^2-R^2}{1+S^2} + \frac{1}{1+S^2} (-v-v^3).$$
Thus, $w_Y$ and $z_X$ can be write as
$$ \displaystyle w_Y = \frac{c'}{4c^2} \frac{R^2 - S^2}{1+R^2} \frac{q}{1+S^2}- \frac{q}{2c} \frac{1}{1+S^2}\frac{1}{1+R^2} (v+v^3),$$
$$ \displaystyle z_X = \frac{c'}{4c^2} \frac{S^2 - R^2}{1+S^2} \frac{p}{1+R^2}- \frac{p}{2c} \frac{1}{1+S^2}\frac{1}{1+R^2} (v+v^3).$$
So
\begin{equation} \label{eq:0219}
\left\{
       \begin{array}{ll}
        \displaystyle   w_Y = \frac{c'}{8c^2} (\cos  z  - \cos  w )q - \frac{q}{8c} (v+v^3)(1+\cos  z )(1+\cos  w ), \\[10pt]
        \displaystyle   z_X = \frac{c'}{8c^2} (\cos  w  - \cos  z )p - \frac{p}{8c} (v+v^3)(1+\cos  z )(1+\cos  w ).
        \end{array}  
        \right.
\end{equation}
By using (\ref{eq:0215}) and (\ref{eq:0218}),
\begin{align*}
\displaystyle p_t - c p_x &= \frac{1}{X_x} 2 R (R_t - c R_x) - \frac{1}{X_x^2} [(X_x)_t - c (X_x)_x] (1+R^2) \\
 \displaystyle   &= \frac{c'}{2c} \frac{p}{1+R^2} [S(1+R^2) - R(1+S^2)] - \frac{p}{1+R^2} R (v+v^3), \\
\displaystyle q_t + c q_x &= \frac{1}{- Y_x} 2 S (S_t - c S_x) - \frac{1}{- Y_x^2} [(-Y_x)_t + c(-Y_x)_x] (1+S^2)\\
\displaystyle   &= \frac{c'}{2c} \frac{q}{1+S^2} [R(1+S^2) - S(1+R^2)] - \frac{q}{1+S^2} S (v + v^3).
\end{align*}
By applying (\ref{eq:0216}),
\begin{align*}
   p_t - c p_x &= -2 c Y_x p_Y ,\\
   q_t + c q_x &= 2 c X_x q_X.
\end{align*}
And thus,
\begin{align*}
p_Y  &= (p_t - c p_x) \frac{1}{-2 c Y_x} = (p_t - c p_x) \frac{1}{2c} \frac{q}{1 + S^2} \\   
     &= \frac{c'}{8c^2} [\sin  z  - \sin  w ] p q - \frac{1}{8c} p q \sin  w  (v + v^3) (1 + \cos  z  ), \\
q_X  &= (q_t + c q_x) \frac{1}{2c X_x} = (q_t + c q_x) \frac{1}{2 c} \frac{p}{1+R^2} \\
     &= \frac{c'}{8c^2} [\sin  w  - \sin  z ] pq - \frac{1}{8c} pq \sin  z  (v+v^3)(1+\cos  w ).
\end{align*}
So, the following identities hold:
\begin{equation} \label{eq:0220}
\left\{
       \begin{array}{ll}
          \displaystyle p_Y = \frac{c'}{8c^2} [\sin  z  - \sin  w ] p q - \frac{1}{8c} p q \sin  w  (v + v^3) (1 + \cos  z ), \\[10pt]
          \displaystyle q_X = \frac{c'}{8c^2} [\sin  w  - \sin  z ] pq - \frac{1}{8c} pq \sin  z  (v+v^3)(1+\cos  w ).
        \end{array}  
        \right.
\end{equation}
Also, we plug in $f=v$ into the equation (\ref{eq:0216}) and get
\begin{equation} \label{eq:0221}
\left\{
       \begin{array}{ll}
           \displaystyle  v_X= (v_t + c v_x) \frac{1}{2c} \frac{p}{1+R^2} = \frac{1}{2c}(\tan \frac{w}{2}\cos ^2\frac{w}{2})p = p\frac{1}{4c} \sin  w  , \\[10pt]
           \displaystyle  v_Y = (v_t - c v_x) \frac{1}{2c} \frac{q}{1+S^2} = \frac{1}{2c}(\tan \frac{z}{2}\cos ^2 \frac{z}{2})q = q\frac{1}{4c} \sin  z  .
        \end{array}  
        \right.
\end{equation}
Combining (\ref{eq:0219}), (\ref{eq:0220}), and (\ref{eq:0221}), we obtain a semi-linear hyperbolic system from the non-linear equation (\ref{eq:011}). This system uses X, Y as independent variables with smooth coefficients for the variables $v,w,z,p,q$
\begin{equation} \label{eq:0222}
\left\{
       \begin{array}{ll}
        \displaystyle   w_Y = \frac{c'}{8c^2} (\cos  z  - \cos  w )q - \frac{q}{8c} (v+v^3)(1+\cos  z )(1+\cos  w ), \\[10pt]
         \displaystyle  z_X = \frac{c'}{8c^2} (\cos  w  - \cos  z )p - \frac{p}{8c} (v+v^3)(1+\cos  z )(1+\cos  w ),
        \end{array}  
        \right.
\end{equation}
\begin{equation} \label{eq:0223}
\left\{
       \begin{array}{ll}
        \displaystyle   p_Y = \frac{c'}{8c^2} [\sin  z  - \sin  w ] p q - \frac{1}{8c} p q \sin  w  (v + v^3) (1 + \cos  z ), \\[10pt]
         \displaystyle  q_X = \frac{c'}{8c^2} [\sin  w  - \sin  z ] pq - \frac{1}{8c} pq \sin  z  (v+v^3)(1+\cos  w ),
        \end{array}  
        \right.
\end{equation}

\begin{equation} \label{eq:0224}
\left\{
       \begin{array}{ll}
           \displaystyle  v_X= \frac{p}{4c} \sin  w  , \\[10pt]
           \displaystyle  v_Y = \frac{q}{4c} \sin  z  ,
        \end{array}  
        \right.
\end{equation}
\indent The system (\ref{eq:0222})-(\ref{eq:0224}) should have non-characteristic boundary conditions related to (\ref{eq:012}). From (\ref{eq:012}), $v_0$ and $v_1$ determine the initial values of $R$ and $S$ at time $t=0$. We denote the curve $\gamma$ as the line in $(X,Y)$ plane at time $t=0$, say
$$ Y = \varphi (X), \,\,\,\,\,\,\,\,\,\, X \in \mathbb{R}.$$
And $ Y = \varphi (X)$ if and only if for some $x \in \mathbb{R},$
$$ X = \int_0^x (1+R^2(0,x)) dx, \,\,\,\,\,\,\,\,\,\,\, Y = \int_x^0 (1 + S^2(0,x)) dx.$$
By the assumptions of the Theorem 1.1, $v_0 \in H^1 \cap L^4, v_1 \in H^1 \cap L^4$. This implies that $R \in H^1 \cap L^4$ and $S \in H^1 \cap L^4$. Moreover, in this case, we let 
\begin{equation} \label{eq:0225}
 \mathcal{E}_0 \coloneqq \frac{1}{4} \int [R^2 (0,x) + S^2(0,x)] dx < \infty .
\end{equation}
Thus, 
\begin{equation} \label{eq:0226}
 X(x) \coloneqq \int_0^x (1+R^2(0,y)) dy, \,\,\,\,\,\,\,\,\,\,\, Y(x) \coloneqq \int_x^0 (1 + S^2(0,y)) dy.
\end{equation}
are absolutely continuous and well defined functions. Further more, by observing (\ref{eq:0226}), $X$ is  increasing and $Y$ is decreasing. So, we conclude that the map $X \mapsto \varphi (X)$ is continuous and decreasing. And from (\ref{eq:0225}),  
$$ | X + \varphi(X) | \leq 4 \mathcal{E}_0.$$
Since $(t,x) \in [ 0, \infty ) \times ( -\infty, \infty )$, so our new independent variables $(X,Y) \in \Omega^+$, and the domain is defined as 
$$ \Omega^+ \coloneqq \{(X,Y): Y \geq \varphi (X)\},$$
along the curve 
$$ \gamma \coloneqq \{ (X,Y): Y = \varphi (X)\}.$$ 
We can have the following boundary data $(\bar{w}, \bar{z}, \bar{p}, \bar{q}, \bar{v}) \in L^{\infty}$,
\begin{equation} \label{eq:0227}
   \left\{
       \begin{array}{ll}
            \bar{w} = 2 \arctan (R(0,x)) ,\\[10pt]
            \bar{z} = 2 \arctan (S(0,x)),
        \end{array}  
        \right.
\end{equation}
\begin{equation} \label{eq:0228}
   \left\{
       \begin{array}{ll}
            \bar{p} \equiv 1 ,\\
            \bar{q} \equiv 1 ,
        \end{array}  
        \right.
\end{equation}
\begin{equation} \label{eq:0229}
         \bar{v} = v_0(x).
\end{equation}
\section{\textbf{Construct the integral solution}} 
We prove the global existence and uniqueness for the semi-linear system (\ref{eq:0222}) - (\ref{eq:0224}) in this section.
\begin{theorem} \label{main}
If the assumptions of Theorem 1.1 holds, then the semi-linear system (\ref{eq:0222}) - (\ref{eq:0224}) with the boundary conditions (\ref{eq:0227}) - (\ref{eq:0229}) has a unique solution for all $(X, Y) \in \mathbb{R} \times \mathbb{R}$.
\end{theorem}
\indent We construct the solution on the region $\Omega^{+}$ which is the case that $Y \geq \varphi (X)$. The proof of the solution on the  $\Omega^{-}$ which is the case that $Y \leq \varphi (X)$ can be construct in the similar way. We show the Lipschitz condition for the system (\ref{eq:0222}) - (\ref{eq:0224}). To make sure the solution is defined in the region $\Omega^+$, we need to construct some priori bounds. So that we can show that $p,q$ are bounded. The Lipschitz condition can be derived as follows. From the energy conservation equations (\ref{eq:015}), (\ref{eq:016}), we denote the following constants:
\begin{equation} \label{eq:0033}
\begin{split}
C_1 &\coloneqq \sup\limits_{v \in \mathbb{R}} \left|\frac{c'(v)}{4c^2(v)} \right| < \infty , \\
  \displaystyle K_1 &\coloneqq \sup\limits_{v \in \mathbb{R}} \int {\frac{v^2}{2} + \frac{v^4}{4}} dx < \infty ,\\
  \displaystyle K_0 &\coloneqq \sup\limits_{v \in \mathbb{R}} v + v^3 < \infty .
\end{split}
\end{equation}
From (\ref{eq:0223}), 
\begin{equation} \label{eq:033}
    \begin{split}
    q_X + p_Y =& \frac{1}{2} \left[ \left(\frac{v^2}{2} + \frac{v^4}{4} \right)q(1+\cos  z )\right]_X + \frac{1}{2} \left[\left(\frac{v^2}{2} + \frac{v^4}{4}\right)p(1+\cos  w )\right]_Y \\ 
   & - \frac{1}{2} {\left(\frac{v^2}{2} + \frac{v^4}{4}\right)\frac{c'}{8c}(\sin  w -\sin  z )(\cos  z  - \cos  w )\left(\frac{1}{c}+1\right)}.
    \end{split}
\end{equation}
\indent We construct a closed curve $\Sigma$ for every $(X,Y) \in \Omega^+$ with the vertical line segment connect $(X,Y)$ with $(X ,\varphi(X))$, the horizontal line segment connect $(X,Y)$ with $(\varphi^{-1}(Y), Y)$, and a part of the boundary $\gamma = {Y= \varphi(X)}$ connecting $(X ,\varphi(X))$ with $(\varphi^{-1}(Y), Y)$. The closed curve $\Sigma = \Gamma_1 + \Gamma_2 + \Gamma_3$. 
\begin{figure}
  \centering
      \includegraphics[scale=0.7]{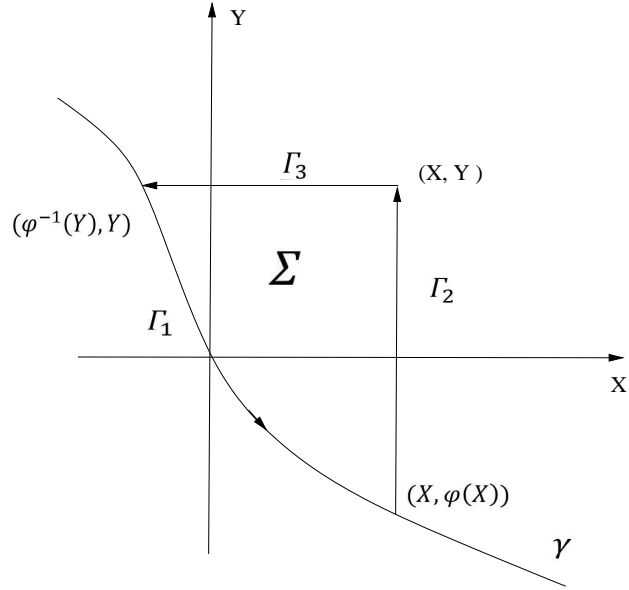}
      \caption{The closed curve $\Sigma$.}
\end{figure}
\indent From (\ref{eq:0033}), we compute  $\displaystyle \iint q_X + p_Y dA = \int -p dX + \int q dY$ and denote that 
\begin{align}
   &\iint q_X + p_Y dA = \int -p dX + \int q dY, \\
   &Q_X \coloneqq \frac{1}{2} \left[\left(\frac{v^2}{2}+\frac{v^4}{4}\right)q(1+\cos  z )\right]_X, \\
   &P_Y \coloneqq \frac{1}{2} \left[\left(\frac{v^2}{2}+\frac{v^4}{4})p(1+\cos  w \right)\right]_Y ,\\
   &\xi \coloneqq \left(\frac{v^2}{2}+\frac{v^4}{4}\right)\frac{c'}{8c}(\sin w - \sin z)(\cos z - \cos w)\left(\frac{1}{c}+1\right).
\end{align}
So
\begin{align*}
  \iint_{\Sigma} p_y + q_X dX dY &= \int_{\Sigma} -p dX + \int_{\Sigma} q dY \\
                            &= \iint_{\Sigma} Q_X + P_Y -\frac{1}{2} \xi dX dY,
 \end{align*}  
By Green's Theorem,
 \begin{align*}         
 \int_{\Sigma} -p dX + \int_{\Sigma} q dY &= \int_{\Sigma} -P dX + \int_{\Sigma} Q dY - \frac{1}{2} \iint_{\Sigma} \xi dX dY. 
\end{align*} 
Thus,  
\begin{align*}
   \frac{1}{2} \iint_{\Sigma} \xi dX dY &= \int_{\Sigma} p - P dX + \int_{\Sigma} Q - q dY, \\
                                         &= \int_{\Sigma} p - \frac{1}{2}\left(\frac{v^2}{2} + \frac{v^4}{2}\right)p(1+\cos w) dX + \int_{\Sigma} -q + \frac{1}{2}\left(\frac{v^2}{2} + \frac{v^4}{2}\right)q(1+\cos z) dX. 
\end{align*}
Since $\Sigma = \Gamma_1 + \Gamma_2 + \Gamma_3$ is a closed curve, so we compute the integral of $\Gamma_1$ directly and in the way $\Gamma_1 = - (\Gamma_2 + \Gamma_3)$.
\begin{equation} \label{eq:039}
\int_{\Gamma_1} 1 - \frac{1}{2} \left(\frac{v^2}{2}+ \frac{v^4}{4}\right)(1+ \cos w)dX + \int_{\Gamma_1} -1 + \frac{1}{2}\left(\frac{v^2}{2}+\frac{v^4}{4}\right)(1+\cos z)dY,
\end{equation}
and from (\ref{eq:0218})
$$ dX = \frac{2}{1+\cos w}dx, \,\,\,\,\,\,\,\,  dY = \frac{2}{1+\cos z}dx.$$ 
Thus $ (1+\cos w)dX = 2 dx$, and  $(1+\cos w)dY = 2 dx.$ \\
So (3.8) becomes
\begin{align*}
\int_{\Gamma_1} 1 - \frac{1}{2} \left(\frac{v^2}{2}+  \frac{v^4}{4}\right)(1+ \cos w)dX + &\int_{\Gamma_1} -1 + \frac{1}{2}\left(\frac{v^2}{2}+\frac{v^4}{4}\right)(1+\cos z)dY \\
                                          & \leq 2 (|X| + |Y| + 4 \mathcal{E}_0) + K_1.
\end{align*}
And also,
\begin{align*}
\int_{\Gamma_2} 1 - \frac{1}{2} \left(\frac{v^2}{2}+ \frac{v^4}{4}\right)(1+ \cos w)dX +  &\int_{\Gamma_2} -1 + \frac{1}{2}\left(\frac{v^2}{2}+\frac{v^4}{4}\right)(1+\cos z)dY \\
& \leq 0 - Y + \varphi(X) + \frac{K_1}{2}, \\
\int_{\Gamma_3} 1 - \frac{1}{2} \left(\frac{v^2}{2}+ \frac{v^4}{4}\right)(1+ \cos w)dX + &\int_{\Gamma_3} -1 + \frac{1}{2}\left(\frac{v^2}{2}+\frac{v^4}{4}\right)(1+\cos z)dY \\
 & \leq \varphi^{-1}(Y) - X  - 0 + \frac{K_1}{2}.
\end{align*}
As a result,
\begin{equation} \label{eq:0310}
\int^X_{\varphi^{-1}(Y)} p (X',Y)dX' + \int^Y_{\varphi(X)} q (X,Y') dY' \leq 2 (|X|+|Y|+4\mathcal{E}_0) + K_1.
\end{equation}
By observing the boundary conditions (\ref{eq:0227}) - (\ref{eq:0229}), $p,q > 0$. And by (\ref{eq:0223}),
$$ p_Y =\frac{1}{8c}pq \left\{ \frac{c'}{c} [\sin z - \sin w] - \sin w (v+v^3) (1+\cos z) \right\}, $$
\begin{align*}
p(X,Y) &= exp \left\{ \int^Y_{\varphi(X)} \frac{1}{8c} \frac{c'}{c} [\sin z - \sin w] - \sin w (v+v^3) (1+\cos z) q(X,Y')dY' \right\} \\
       &\leq exp \left\{ C_1 \int^Y_{\varphi(X)} q(X,Y') dY' \right\} \\
       &\leq exp \left\{ ~ 2 C_1 (|X|+|Y|+ 4 \mathcal{E}_0) + C_1 K_1 \right\}.
\end{align*}
Similarly, we have 
$$ q(X,Y) \leq  exp \{ ~2 C_1 (|X|+|Y|+ 4 \mathcal{E}_0) + C_1 K_1 \}.$$ 
 \indent Now, we show that on any bounded sets in $X$-$Y$ plane, we can construct the solution for the system of the equations (\ref{eq:0222}) - (\ref{eq:0224}) with boundary condition (\ref{eq:0227}) -(\ref{eq:0229}) by the fixed point of a constructive map. For any $r > 0$, we can construct a bounded domain
 $$ \Omega_r \coloneqq \{(X,Y): Y \leq \varphi (X), X \leq r, Y \leq r\}.$$
 And also introduce the function space :
 \begin{equation} \label{eq:0311}
 \displaystyle
  \Lambda_r \coloneqq \{f: \Omega_r \mapsto \mathbb{R}: \| f \|_{*} \coloneqq \,\, \text{ess} \hspace{-4mm} \sup\limits_{(X,Y) \in \Omega_r} e^{-K(X+Y)} |f(X,Y) | < \infty \}.
 \end{equation}
Where $K$ is a suitably big constant and it will be determined later. And for $(w,z,p,q,v) \in \Lambda_r$, we construct a map $ \tau(w,z,p,q,v) = (\tilde{w},\tilde{z},\tilde{p},\tilde{q},\tilde{v})$. And this map is define as follows. 
  \begin{equation} \label{eq:0312}
\left\{
       \begin{array}{ll}
          \tilde{w}(X,Y) =  \bar{w}(X, \varphi (X)) + \int^Y_{\varphi(X)} \frac{c'}{8c^2} (\cos  z  - \cos  w )q - \frac{q}{8c} (v+v^3)(1+\cos  z )(1+\cos  w ) dY ,\\[10pt]
          \tilde{z} (X,Y) = \bar{z} (\varphi^{-1} (Y), Y) + \int^X_{\varphi^{-1} (Y)} \frac{c'}{8c^2} (\cos  w  - \cos  z )p - \frac{p}{8c} (v+v^3)(1+\cos  z )(1+\cos  w ) dX,
        \end{array}  
        \right.
\end{equation}
\begin{equation} \label{eq:0313}
\left\{
       \begin{array}{ll}
          \tilde{p} (X,Y) = 1 + \int^Y_{\varphi(X)} \frac{1}{8c}pq \left\{ \frac{c'}{c} [\sin z - \sin w] - \sin w(v+v^3)(1 + \cos z) \right\} dY,\\[10pt]
          \tilde{q} (X,Y) = 1 + \int^X_{\varphi^{-1}(Y)} \frac{1}{8c}pq \left\{ \frac{c'}{c} [\sin w - \sin z] - \sin z(v+v^3)(1 + \cos w) \right\} dX,
        \end{array}  
        \right.
\end{equation}
\begin{equation} \label{eq:0314}
            \tilde{v} (X,Y)=  \bar{v} (X, \varphi (X))+ \int^Y_{\varphi (X)} \frac{1}{4c} \sin  z  q dY .
\end{equation}
We want to prove the uniform Lipschitz condition. First, we define 
$$ \Phi_r \coloneqq \Lambda_r \times \Lambda_r \times \Lambda_r \times \Lambda_r \times \Lambda_r .$$
For some properly chosen distance $D : \Phi_r \times \Phi_r \mapsto \mathbb{R}$, we want to show that 
$$ D ( (\tilde{w_1},\tilde{z_1},\tilde{p_1},\tilde{q_1},\tilde{v_1}), (\tilde{w_2},\tilde{z_2},\tilde{p_2},\tilde{q_2},\tilde{v_2})) < L \times D((w_1,z_1,p_1,q_1,v_1) , (w_2,z_2,p_2,q_2,v_2)).$$
The Lipschitz constant $L$ satisfies $L \leq 1$. In fact, we define the distance as: 
\begin{align*}
 D( (\tilde{w_1},\tilde{z_1},\tilde{p_1},\tilde{q_1},\tilde{v_1}),  (\tilde{w_2},\tilde{z_2},\tilde{p_2},\tilde{q_2},\tilde{v_2})) \coloneqq \max &\{ \|\tilde{w_1} - \tilde{w_2}  \|_*\,\, ,\,\, \|\tilde{z_1}
  - \tilde{z_2}  \|_*\,\, , \\
  & \,\, \|\tilde{p_1} - \tilde{p_2}  \|_* \,\, , \,\, \|\tilde{q_1} - \tilde{q_2} \|_* \,\, , \,\, \|\tilde{v_1} - \tilde{v_2}  \|_* \},
\end{align*}
and the norm $\| \cdot \|_* $ is defined in (\ref{eq:0311}). \\
\indent A straightforward computation shows that $ \displaystyle{L = \frac{C(\mathcal{E}_0,\mathcal{K})}{K}}$, where $C(\mathcal{E}_0,\mathcal{K})$ is a constant depends on $\mathcal{E}_0$ and $\mathcal{K}$. By choosing $K$ sufficiently large, we can guarantee $L < 1$. Hence, the uniform Lipschitz condition is proved. By the fixed point theorem, the solution in the $X$-$Y$ plane exists and is unique. \qed \\
 \indent If the initial data in (1.2) are smooth, then the solutions of (\ref{eq:0222}) - (\ref{eq:0224}) with boundary condition (\ref{eq:0227}) - (\ref{eq:0229}) are smooth functions with variables $(X,Y)$. Also, if there is a sequence of smooth functions $(v^m_0(x),v^m_1(x))_{m\geq 1}$ with the following conditions: 
      \begin{center}
       $ v^m_0(x) \to v_0(x) $ ,   $ v^m_1(x) \to v_1(x) $,    $(v^m_0(x))_x \to (v_0(x))_x $,        
      \end{center}
      uniformly on a compact subset of $\mathbb{R}$. Then
      $$ (p^m,q^m,w^m,z^m,v^m) \to (p,q,w,z,v),$$
      uniformly on some bounded subsets of $X$-$Y$ plane.
\section{\textbf{Weak solutions }} 
\indent In this section, we construct a map $ v(X,Y) \to v(t,x)$. That is to write $(X,Y)$ in terms of $(t,x)$ so we obtain a solution to the Cauchy problem (\ref{eq:011}), (\ref{eq:012}). The map $(X,Y) \mapsto (t,x)$ can be obtain in the following way. We plug in $f = x$ and $ f = t$ into the equation (\ref{eq:0216}), and get
\begin{equation} \label{eq:041}
   \left\{
       \begin{array}{ll}
               c = 2 c X_x x_X ,\\
               -c = -2 c Y_x x_Y ,\\
               1 = 2 c X_x t_X ,\\
               1 = -2 c Y_x t_Y.
        \end{array}  
        \right.
\end{equation}
And by applying (\ref{eq:0218}) we have
\begin{equation} \label{eq:042}
   \left\{
       \begin{array}{ll}
          \displaystyle X_x = \frac{2}{(1+\cos w)p},\\[10pt]
          \displaystyle Y_x = \frac{-2}{(1+\cos z)q}, \\[10pt]
          \displaystyle X_t = \frac{2c}{(1+\cos w)p},\\[10pt]
          \displaystyle Y_t = \frac{-2c}{(1+\cos z)q}. \\
        \end{array}  
        \right.
\end{equation}
 \indent We assume that the partial derivatives above valid for points that $w,z \neq - \pi$. Thus, we have
\begin{equation} \label{eq:043}
   \left\{
       \begin{array}{ll}
         \displaystyle   x_X = \frac{1}{2 X_x} = \frac{(1+\cos w)p}{4}, \\[10pt]
          \displaystyle  x_Y = \frac{1}{2 Y_x} = \frac{-(1+\cos z)q}{4},
        \end{array}  
        \right.
\end{equation}
\begin{equation} \label{eq:044}
   \left\{
       \begin{array}{ll}
         \displaystyle   t_X = \frac{1}{2c X_x} = \frac{(1+\cos w)p}{4c} ,\\[10pt]
          \displaystyle  t_Y = \frac{1}{-2c Y_x} = \frac{(1+\cos z)q}{4c}.
        \end{array}  
        \right.
\end{equation}
A computation shows that $x_{XY} = x_{YX}$ and $t_{XY} = t_{YX}$
\begin{align*}
    x_{XY}  =& \frac{(1+\cos w)p_Y}{4} - \frac{p \sin w w_Y}{4} \\
             = &\frac{c' pq}{32c}[\sin z - \sin w + \sin (z-w)] , \\
   x_{YX}   = &\frac{(1+\cos z)q_X}{4} - \frac{q \sin z z_Y}{4} \\
             = &\frac{c' pq}{32c}[\sin z - \sin w + \sin (z-w)] .
\end{align*}
So, $ x_{XY} = x_{YX}$. \\
\indent And similarly, we can compute that $t_{XY} = t_{YX}$. Thus, the two equation in (\ref{eq:043}) are equivalent: $x_{XY} = x_{YX}$. And the two equation in (\ref{eq:044}) are equivalent since $t_{XY} = t_{YX}$.  We can recover the solution in terms of $(t,x)$ with function $x = x(X,Y)$ by integrating one of the equation in (\ref{eq:043}). Also, we can write $t = t(X,Y)$ by integrating one of the equation in (\ref{eq:044}).  \\
\indent Next, we prove that the function $v$ is a weak solution to (\ref{eq:011}). From (\ref{eq:014}), we want to show that 
\begin{align*}
0  &=  \iint \phi_t v_t - [c(v) \phi]_x [c(v) v_x] - \frac{\phi}{2} (v+v^3) dx dt.
\end{align*}
In fact, it is equivalent to prove:
\begin{align*}
0   =& \iint (v_t + c v_x) [\phi_t - (c(v)\phi)_x] + (v_t - c v_x) [\phi_t +(c(v)\phi)_x] - \phi(v + v^3)dx dt \\ 
   =&  \iint -(\frac{\sin w}{2} p)_Y \phi - (\frac{\sin z}{2} q)_X \phi + \frac{c' pq}{8c} [\sin w \frac{1+\cos z}{2} - \sin z \frac{1+\cos w}{2}]\phi (\tan  \frac{z}{2} - \tan  \frac{w}{2}) dX dY \\
   & - \iint \phi (v + v^3) dx dt \\
   \coloneqq & \text{I} + \text{II}.  
\end{align*}
We define I and II later, and where in the last step, we have used (\ref{eq:042}),
\begin{align*}
dx dt = \displaystyle 
\begin{vmatrix}
\displaystyle \frac{dx}{dX} & \displaystyle \frac{dx}{dY} \\[10pt]
\displaystyle \frac{dt}{dX} & \displaystyle \frac{dt}{dY}
\end{vmatrix}
dX dY = \displaystyle \frac{pq}{2c(1+R^2)(1+S^2)} dX dY.
\end{align*}
And used  the following identities derived from (\ref{eq:0224}), 
\begin{equation} \label{eq:047}
   \left\{
       \begin{array}{ll}
         \displaystyle \frac{1}{1+R^2} = \frac{1 +\cos w}{2} ,\\[10pt]
          \displaystyle \frac{1}{1+S^2} = \frac{1 + \cos z}{2},
        \end{array}  
        \right.
\end{equation}
\begin{equation} \label{eq:048}
   \left\{
       \begin{array}{ll}
         \displaystyle  \frac{R}{1+R^2} = \frac{\sin w}{2} ,\\[10pt]
         \displaystyle  \frac{S}{1+S^2} = \frac{\sin z}{2}.
        \end{array}  
        \right.
\end{equation}
We denote I and II as follows 
\begin{equation} \label{eq:045}
\begin{multlined}
\text{I} =  \iint -\left(\frac{\sin w}{2} p\right)_Y \phi - \left(\frac{\sin z}{2} q\right)_X \phi \\
+ \frac{c' pq}{8c} \left[\sin w \frac{1+\cos z}{2} - \sin z \frac{1+\cos w}{2}\right]\phi \left(\tan  \frac{z}{2} - \tan  \frac{w}{2}\right) dX dY ,
   \end{multlined}
\end{equation}
and 
\begin{equation} \label{eq:046}
\text{II} = \iint \phi (v + v^3) dx dt.
\end{equation}
A computation on I with (\ref{eq:0222}) - (\ref{eq:0224}) shows that
\begin{align*}
 \text{I}  =& \iint -\left(\frac{\cos w}{2} w_Y p + \frac{\sin w}{2} p_Y \right) \phi -\left(\frac{\cos z}{2} z_X q + \frac{\sin z}{2} q_X\right) \phi + \frac{c' p q}{8 c^2} [\cos (w+z) - 1]\phi dX dY \\
 =& \iint  \frac{p q}{16c} \phi (v+v^3) (\cos w + \cos z + 2 + 2 \cos z \cos w  \\
 &+\cos ^2 w \cos z + \cos w \cos ^2 z + \sin ^2 w \cos z + \sin ^2 z \cos  w) dX dY \\
 =& \iint \frac{p q }{8 c} \phi (v+v^3) (1 + \cos z + \cos w + \cos z \cos w) dX dY .
\end{align*} 
A computation on II shows that 
 \begin{align*}
   \text{II} &= \iint \phi (v + v^3) \frac{p q}{ 2 c (1 + s^2)(1 + R^2)} dX dY \\
         &= \iint \frac{p q}{8 c} \phi (v + v^3) (1 + \cos z + \cos w + \cos z \cos w) dX dY.
 \end{align*}
Clearly, $ \text{I}  = \text{II}$. Thus the integral (\ref{eq:014}) holds since 
 \begin{align*}
0  =&  \iint -\left(\frac{\sin w}{2} p\right)_Y \phi - \left(\frac{\sin z}{2} q\right)_X \phi\\
 & + \frac{c' pq}{8c} \left[\sin w \frac{1+\cos z}{2} - \sin z \frac{1+\cos w}{2}\right]\phi \left(\tan  \frac{z}{2} - \tan  \frac{w}{2}\right) dX dY   - \iint \phi (v + v^3) dx dt \\
0  =& ~ \text{I}  - \text{II} ,
\end{align*}
where I is defined in (\ref{eq:045}), and II is defined in (\ref{eq:046}). \\
\indent Next, we define $v$ as a function in terms of the original variables $(t,x)$. We invert the map $(X,Y) \mapsto (t,x)$ and then we have $v(t,x) = v(X(t,x), Y(t,x))$. Given arbitrary $(t^*, x^*)$ in the $t\text{-}x$ plane, we choose arbitrary point $(X^*, Y^*)$ in $X$-$Y$ plane such that $t^* = t(X^*, Y^*)$ and $x^* = x(X^*,Y^*)$. We define that $v(t^*,x^*) = v(X^*,Y^*)$ and assume that there are two different points $(t(X_1, Y_1),x(X_1, Y_1) ) = (t(X_2,Y_2), x(X_2,Y_2)) = (t^*,x^*)$. We consider two cases: case 1: $X_1 \leq X_2, Y_1 \leq Y_2$, and case 2: $X_1 \leq X_2, Y_1 \geq Y_2$. 
\begin{figure}
  \centering
       \includegraphics[scale=0.6]{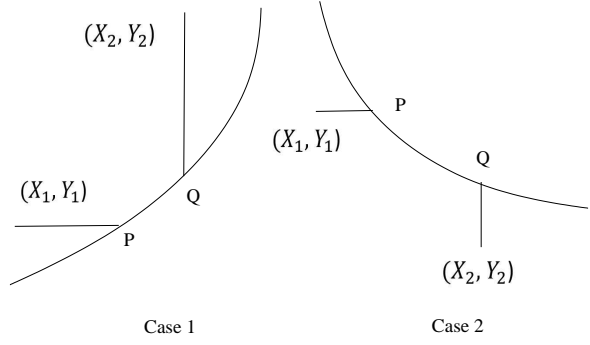}
      \caption{The paths of integration.}
\end{figure}
\indent Case 1: $X_1 \leq X_2, Y_1 \leq Y_2$. We consider the set
  $$  \Gamma_{x^*} \coloneqq \{(X,Y) : x(X,Y) \leq x^*\}.$$
We denote $\partial \Gamma_{x^*}$ as the boundary of $\Gamma_{x^*} $. By (\ref{eq:043}), we observe that $x$ is increasing with $X$ increasing and $x$ is decreasing with $Y$ increasing. Thus, this boundary can be write as a Lipschitz continuous function denoted as $X -Y = \phi (X- Y)$. 
  \indent We construct the Lipschitz continuous curve $\gamma$ with the following properties:
  \begin{itemize}
   \item   a horizontal line segment connecting $(X_1, Y_1)$ with a point $P = (X_P, Y_P)\in \partial \Gamma_{x^*}$  and $Y_P = Y_1$.
   \item  a vertical line segment connecting $(X_2, Y_2)$ with a  point $Q = (X_Q, Y_Q)\in \partial \Gamma_{x^*}$  and $X_Q = X_2$.
   \item  a part of $ \partial \Gamma_{x^*}$.  
  \end{itemize}
Thus, we obtain a Lipschitz continuous parametrization of the curve $\gamma$ : $[\xi_1 , \xi_2] \mapsto \mathbb{R} \times \mathbb{R}$ where the parameter $\xi = X + Y$. By observing, the map $(X,Y) \mapsto (t,x) $ is constant along the curve $\gamma$. And (\ref{eq:043}) - (\ref{eq:044}) implies that
  \begin{equation} \label{eq:049}
  (1+\cos w)X_{\xi} = (1+\cos z)Y_{\xi} = 0, 
  \end{equation}
  From (\ref{eq:049}),
  \begin{equation}  \label{eq:0410}
   \sin w X_{\xi} = \sin z Y_{\xi} = 0.  
  \end{equation}
Thus, by (\ref{eq:0410})
\begin{align*}
  v(X_2,Y_2) - v(X_1, Y_1)  &= \int_{\gamma} (v_X dX + v_Y dY) \\
                            &= \int^{\xi_2}_{\xi_1} (\frac{p\sin w}{4c} X_{\xi} - \frac{q\sin z}{4c}Y_{\xi})d\xi = 0.
\end{align*}
  So our claim for case 1 is proved.\\
  \indent Case 2: $X_1 \leq X_2, Y_1 \geq Y_2$. We consider the set:
  $$\Gamma_{t^*} \coloneqq \{(X,Y) : t(X,Y) \leq t^*\}.$$
  And we do the same process as we did in case 1. Construct $\gamma$ connecting $(X_1, X_2)$ and $(X_2,Y_2)$ as Figure 3 case 2 indicates. \\
\indent Next, we prove the function $v(t,x) = v(X(t,x),Y(t,x))$ is H\"older-$\frac{1}{2}$ continuous on the bounded sets. To prove this, we need to consider characteristic curve such that $t \mapsto x^+(t)$ with $\bar{x}^+ = c(v)$. For some fixed $\bar{Y}$, this can be parametrized by the function $X \mapsto (t(X,\bar{Y}), x(X, \bar{Y}))$. By (\ref{eq:0214}), (\ref{eq:0216}), (\ref{eq:0218}) and (\ref{eq:0224}), 
  \begin{align*}
  \int^{\tau}_0 [v_t + c(v) v_x]^2 dt &= \int^{X_{\tau}}_{X_0} (2 c X_x v_X)^2 \frac{1}{2 X_t} dX\\
                                      &= \int^{X_{\tau}}_{X_0} \frac{p}{2c} \sin ^2(\frac{w}{2}) dX \leq \int^{X_{\tau}}_{X_0} \frac{p}{2c} dX \leq C_{\tau}.
  \end{align*}
  Thus, we obtain that
  \begin{equation} \label{eq:0411}
  \int^{\tau}_0 [v_t + c(v) v_x]^2 dt \leq C_{\tau}.
  \end{equation}
Similarly, we integrate along backward characteristics curves $t \mapsto x^{-} (t)$ and find out that 
  \begin{equation} \label{eq:0412}
  \int^{\tau}_0 [v_t - c(v) v_x]^2 dt \leq C_{\tau} .
  \end{equation}
Thus, since the speed of the characteristic curve is $+ c(v)$ or $- c(v)$ and $c(v)$ is uniformly positive bounded. With the bounds (\ref{eq:0411}) and (\ref{eq:0412}), the function $v(t,x)$ is H\"older-$ \frac{1}{2}$ continuous.  \qed
\section{\textbf{Conserved quantities }} 
\indent This section provides a proof of Theorem 1.2. Recalling (\ref{eq:026}) and (\ref{eq:027}), a straightforward computation shows that 
\begin{align*}
E_t &= \left(\frac{1}{2} v_t^2 + \frac{1}{2} c^2 v_x^2 + \frac{v^2}{4} + \frac{v^4}{8}\right)_t \\
    &= v_{tt} v_t + c c' v_t v_x^2 + c^2 v_x v_{xt} + \frac{1}{2}(v + v^ 3) v_t ,\\
(c^2 M) _x &= (- c^2 v_t v_x)_x = - 2 c c' v_x^2 v_t - c^2 v_{tx} v_x - c^2 v_t v_{xx},\\
E_t + (c^2 M)_x  &= v_t (v_{tt} - c c' v_x - c^2 v_{xx} + \frac{1}{2} (v+ v^3)) = 0,
\end{align*}
and 
\begin{align*}
M_t  &= - v_{tt} v_x - v_t v_{xt}, \\
E_x  &= v_t v_{tx} + c c' v_x v_x^2 + c^2 v_x v_{xx} + \frac{1}{2} (v + v^3) v_x ,\\
M_t + E_x  &= -v_x \left(v_{tt} - c c' v^2_x - c^2 v_{xx} - \frac{1}{2} (v+ v^3)\right) = 0.
\end{align*}
Thus,
\begin{equation} \label{eq:051}
   \left\{
       \begin{array}{ll}
          E_t + (c^2 M)_x = 0 ,\\
          M_t + E_x = 0.
        \end{array}  
        \right.
\end{equation}
Also,
\begin{equation} \label{eq:052}
 dx = \frac{(1+\cos w)p}{4} dX - \frac{(1+\cos z)q}{4} dY, 
\end{equation}
\begin{equation} \label{eq:053}
 dt = \frac{(1+\cos w)p}{4c} dX + \frac{(1+\cos z)q}{4c} dY ,
\end{equation}
which is closed.
We want to show that $Edx - (c^2 M) dt , \,\,\, M dx - E dt$ are closed. Recalling (\ref{eq:0222}) - (\ref{eq:0224}), we write $Edx - (c^2 M) dt , \,\,\, M dx - E dt$ in terms of $X,Y$, and show that they are closed.
\begin{equation} \label{eq:055}
\begin{split}
& E dx - (c^2 M) dt = \\
 &   \left[\frac{(1- \cos w)}{8} + \frac{(1+\cos w)}{32}(2v^2 + v^4)\right]p dX - \left[\frac{(1-\cos z)}{8} + \frac{1+\cos z}{32}(2v^2 + v^4)\right]q dY ,
  \end{split}
\end{equation}
\begin{align*}
&M dx + E dt = \\
& \left\{  \frac{(1 - \cos w)}{8c} + \frac{(1+\cos w)}{32c} (2 v^2 + v^4)  \right\} p dX + \left\{  \frac{1 - \cos z}{8c} + \frac{1 + \cos z}{32c} (2 v^2 + v^4) \right\} q dY .
\end{align*}
And we compute that 
\begin{align*}
&\left\{\left[\frac{(1- \cos w)}{8} + \frac{(1+\cos w)}{32}(2v^2 + v^4)\right]p \right\}_Y \\
= & \frac{-\sin w (2 v^2 + v^4)p}{32} \frac{c'}{8c^2} (\cos z - \cos w) q + \frac{(1+\cos w)(2v + v^3)p}{32} \frac{1}{4c} \sin z q\\
&+ \frac{1 + \cos w}{32} (2 v^2 v^4) \frac{c'}{8c^2} [\sin z - \sin w]pq\\
=& - \left\{ \left[\frac{(1-\cos z)}{8} + \frac{1+\cos z}{32}(2v^2 + v^4)\right]q\right\}_X,
\end{align*}
and 
\begin{align*}
& \left\{ \left[  \frac{(1 - \cos w)}{8c} + \frac{(1+\cos w)}{32c} (2 v^2 + v^4)  \right] p  \right\}_Y  \\
=& ~ \frac{\sin w p}{64c^3} c' (\cos z - \cos w)q - \frac{\sin w p q}{64c^2} (v+v^3)(1+\cos z)(1+\cos w) \\
& + \frac{(1 -\cos w)c'}{64 c^3} [\sin z - \sin w] pq - \frac{(1-\cos w) pq}{64c^2} \sin w(v+v^3)(1+\cos z) \\
& - \frac{\sin w}{32c} p (2v^2 + v^4) \frac{c'}{8 c^2} (\cos z - \cos w)q + \frac{(1+\cos w)}{32c} (2v^2 + v^4)\\ & \frac{c'}{8c^2} [\sin z - \sin w] p q + \frac{1 + \cos w}{ 32c} p (4 v + 4 v^3) \frac{1}{4c} \sin z q \\
=&  \left\{ \left[  \frac{1 - \cos z}{8c} + \frac{1 + \cos z}{32c} (2 v^2 + v^4) \right] q \right\}_X.
\end{align*}
Thus $ \{Edx - (c^2 M) dt \}, \,\,\,\{ M dx - E dt\}$ are closed.
\begin{figure}
  \centering
      \includegraphics[scale=0.65]{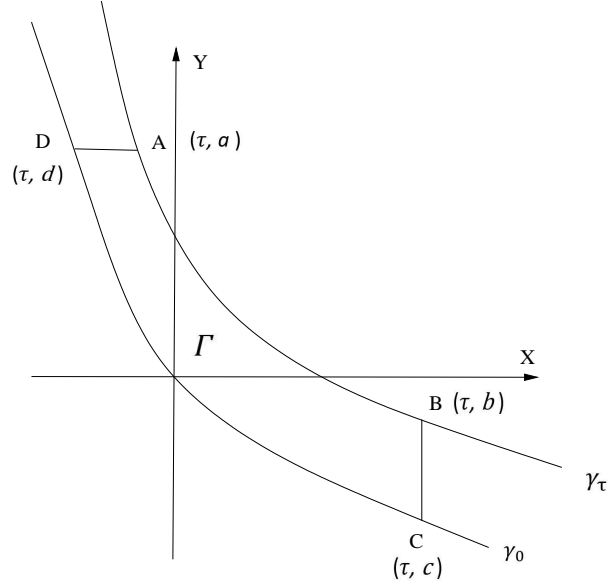}
      \caption{The region $\Gamma$.}
\end{figure} \\
\indent To prove the inequality (\ref{eq:017}), We fixed some $\tau > 0$, and the case $\tau < 0$ is identical. 
We assume that for an arbitrary large $r > 0$. We define the set 
\begin{equation} \label{eq:056}
\Gamma \coloneqq \{(X,Y) : 0 \leq t(X,Y) \leq \tau, X \leq r, Y \leq r\}.
\end{equation}
We form the map $(X,Y) \mapsto (t,x)$ in the following pattern:
$$ A \mapsto (\tau, a), \,\,\,\,\, B \mapsto (\tau, b) , \,\,\,\, C \mapsto (0, c) , \,\,\,\,\ D \mapsto(0, d),$$ 
such that $a < b$ and $c > d$. Then, we can integrate the (\ref{eq:055}) along $ \partial \Gamma$, the boundary of $\Gamma$.
\begin{align*}
&\int_{AB} \left\{ \frac{(1- \cos w)p}{8} + \frac{(1+\cos w)p}{32} (2 v^2 + v^4) \right\} dX - \left\{ \frac{(1-\cos z)q}{8} + \frac{(1+\cos z)q}{32} (2 v^2 + v^4)  \right\} dY \\
 =& \int_{DC} \left\{ \frac{(1- \cos w)p}{8} + \frac{(1+\cos w)p}{32} (2 v^2 + v^4) \right\} dX - \left\{ \frac{(1-\cos z)q}{8} + \frac{(1+\cos z)q}{32} (2 v^2 + v^4)  \right\} dY \\
& - \int_{DA} \left\{ \frac{(1- \cos w)p}{8} + \frac{(1+\cos w)p}{32} (2 v^2 + v^4) \right\} dX \\
& - \int_{CB} \left\{ \frac{(1-\cos z)q}{8} + \frac{(1+\cos z)q}{32} (2 v^2 + v^4)  \right\} dY \\
 \leq & \int_{DC} \left\{ \frac{(1- \cos w)p}{8} + \frac{(1+\cos w)p}{32} (2 v^2 + v^4) \right\} dX - \left\{ \frac{(1-\cos z)q}{8} + \frac{(1+\cos z)q}{32} (2 v^2 + v^4)  \right\} dY \\
 \leq& \int^c_d \frac{1}{2} \left[v^2_t (0, x) + c^2 (v(0,x)) v^2_x (0,x) + \frac{1}{2} (v(0,x) + v^3(0,x)) \right] dx.
\end{align*}
Also,
\begin{align*}
& \int^b_a \frac{1}{2} \left[v^2_t (0, x) + c^2 (v(0,x)) v^2_x (0,x) + \frac{1}{2} (v(0,x) + v^3(0,x))\right] dx\\ \nonumber
 =& \int_{AB \cap \{ \cos w \neq - 1\}} \left\{ \frac{(1- \cos w)p}{8} + \frac{(1+\cos w)p}{32} (2 v^2 + v^4) \right\} dX \\
& - \left\{ \frac{(1-\cos z)q}{8} + \frac{(1+\cos z)q}{32} (2 v^2 + v^4)  \right\} dY \\
 \leq& ~ \mathcal{E}_0.
\end{align*}
Let $r \to \infty$, $a \to - \infty$ and $b \to + \infty$. We conclude that $\mathcal{E}(t) \leq \mathcal{E}_0$. Thus, the inequity  (\ref{eq:017}) is proved.
\begin{figure}
  \centering
      \includegraphics[scale=0.7]{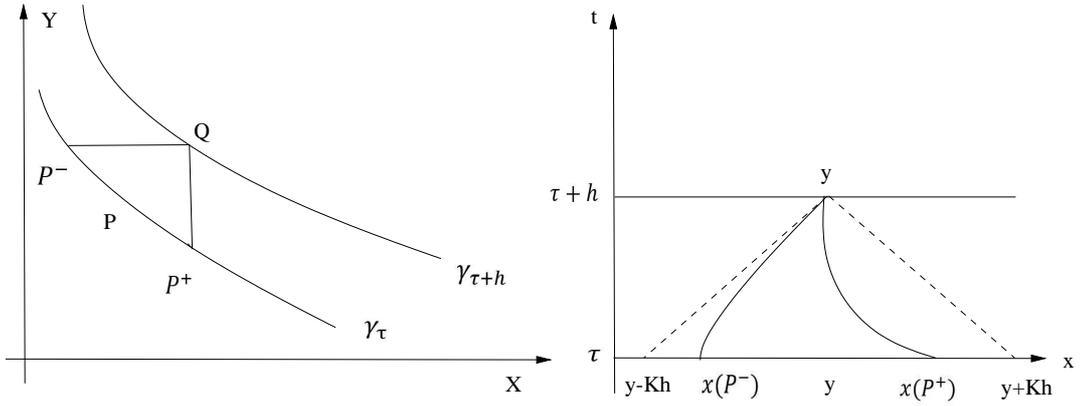}
      \caption{Proving Lipschitz condition.}
\end{figure} \\
\indent Now, we prove the Lipschitz condition on the map $t \mapsto v (t , \cdot)$ in the $L^2$ distance. First, for any fixed time $\tau$, we define $\mu_{\tau} \coloneqq \mu_{\tau}^- + \mu_{\tau}^+$ and $\mu_{\tau }$ is the positive measure on the real lines. We define $\mu_{\tau}^-, \mu_{\tau}^+$ as follows. \\
we define $\Gamma_{\tau} \coloneqq \{  (X,Y): t(X,Y) \leq \tau \}$ and let $\gamma_{\tau}$ be the boundary of $\Gamma_{\tau}$. \\
\indent For any open interval $] a,b[$, we define $A = (X_A, Y_A), B= (X_B, Y_B)$ be points on the $\gamma_{\tau}$ such that \\
\indent $x(A) = a$, and $X_P - Y_P \leq X_A - Y_A$ for all points $P \in \gamma_{\tau}$ and $x(P) \leq a$,\\
\indent $x(B) = b$, and $X_B - Y_B \leq X_P - Y_P$ for all points $P \in \gamma_{\tau}$ and $x(P) \geq b$.\\
Then we have 
\begin{equation} \label{eq:057}
 \mu_{\tau} \coloneqq \mu^-_{\tau}(]a,b[) + \mu^+_{\tau}(]a,b[),
\end{equation}
and define that in the general case
\begin{equation} \label{eq:058}
 \mu^-_{\tau}(]a,b[)   \coloneqq \int_{AB} \left\{ \frac{(1- \cos w)p}{8} + \frac{(1+\cos w)p}{32} (2 v^2 + v^4) \right\} dX ,
\end{equation}
\begin{equation} \label{eq:059}
 \mu^+_{\tau} (]a,b[)  \coloneqq \int_{AB} - \left\{ \frac{(1-\cos z)q}{8} + \frac{(1+\cos z)q}{32} (2 v^2 + v^4) \right\} dY.
\end{equation}
In the smooth case:
\begin{equation} \label{eq:0510}
\mu^-_{\tau} (] a, b [) \coloneqq \frac{1}{4} \int^b_{a} R^2 (\tau, x) dx ,
\end{equation}
 \begin{equation} \label{eq:0511}
 \mu^+_{\tau} (] a, b [) \coloneqq \frac{1}{4} \int^b_a S^2 ( \tau , x) dx.
 \end{equation}
Clearly, $\mu^+$ and $\mu^-$ are bounded positive measure. For all $\tau$, we have $\mu_{\tau} (\mathbb{R}) = \mathcal{E}_0$ by (\ref{eq:055}). By (\ref{eq:0510})- (\ref{eq:0511}) and (\ref{eq:027}) we compute that 
\begin{align*}
\int^b_a c^2 u_x^2 dx= \int^b_a \frac{c^2(R - S)^2}{4c^2} dx &= \int^b_a \frac{R^2 - 2RS + S^2}{4} dx \\
                                                 &\leq \int^b_a \frac{R^2 + S^2}{2} dx = 2 \mu(]a,b[).
\end{align*}
Thus, for arbitrary $a,b$ with $a < b$,
\begin{equation} \label{eq:0512}
|v(\tau,b) - v(\tau,a)|^2 \leq |b-a| \int^b_a v_x(\tau,y) dy \leq |b-a|2\mathcal{K}^2 \mu_{\tau} (]a,b[).
\end{equation}\\
\indent For given $y \in R$ and $h > 0$, our goal is to estimate the $|v(\tau + h,y) - v(\tau, y)|$. We first denote that $\Gamma_{\tau+h}$ as the set $\Gamma_{\tau+h} \coloneqq \{(X,Y): t(X,Y ) \leq \tau + h\}$ and denote that $\gamma_{\tau+h}$ to be the boundary of the set $\Gamma_{\tau+h}$. \\
\indent Let $P=(P_X, P_Y)$ be points on $\gamma_{\tau+h}$ (as the figure 5(a) shows) such that $x(P) = y$, and $X_{\tilde{P}} - Y_{\tilde{P}} \leq X_P - Y_P$ for all $\tilde{P} \in \gamma_{\tau}$, $x(\tilde{P}) \leq x(P)$.\\
\indent Let $Q = (Q_X, Q_Y)$ be points on $\gamma_{\tau + h}$ such that $x(Q) = y$ and $X_{\tilde{Q}} - Y_{\tilde{Q}} \leq X_Q - Y_Q$ for all $\tilde{Q} \in \gamma_{\tau + h}$, $x(\tilde{Q}) \leq x(Q)$. \\
\indent So $X_P \leq X_Q$ and $Y_P \leq Y_Q$. Let $P^+ = (X_Q, Y^+) \in \gamma_{\tau}$ and $P^- = (X^-, Y_Q) \in \gamma_{\tau}$.\\
\indent As shown in the figure 5, since the point $(\tau, x(P^+))$ lies on some characteristic curve with the speed $c(v) \leq \mathcal{K}$ and go through the point $(\tau + h, y)$, so $x(P^+) \in ] y, y + \mathcal{K}h[$.
Also, $x(P^-) \in ]  y - \mathcal{K}h, y[$, since point $(\tau, y)$ lies on some characteristic curve with the speed $- c(v) \geq -\mathcal{K}$ and go through the point $(\tau + h , y)$.\\
Thus, and by (\ref{eq:0224}), we can compute that 
\begin{align*}
 | v(Q) - v (P^+)| &\leq \int^{Y_Q}_{Y^+} |v_Y (X_Q, Y)|dY \\
                   & = \int^{Y_Q}_{Y^+} \left(\frac{1 + \cos z}{4c} q\right)^{\frac{1}{2}}  \left(\frac{1 - \cos z}{4c} q\right)^{\frac{1}{2}} dY\\
                   &\leq \left(\int^{Y_Q}_{Y^+} \frac{1 + \cos z}{4c} q dY\right)^{\frac{1}{2}} \left(\int^{Y_Q}_{Y^+} \frac{1 - \cos z}{4c} q dY\right)^{\frac{1}{2}}   \\
                   &\leq \left(\int^{Y_Q}_{Y^+} \frac{1 + \cos z}{4c} q dY + \frac{1 + \cos w}{4c} p dX \right)^{\frac{1}{2}} \left(\int^{Y_Q}_{Y^+} \frac{1 - \cos z}{4c} q dY + \frac{1 - \cos w}{4c} p dY\right)^{\frac{1}{2}}\\
                   & \leq \left(\int^{\tau+h}_{\tau} 1 dt \right)^{\frac{1}{2}} \left(\int^{P^+}_{P^-} \frac{1 - \cos z}{4c} q dY + \frac{1 - \cos w}{4c} p dY\right)^{\frac{1}{2}} \\
                  & \leq h^{\frac{1}{2}} \left(\int^{P^+}_{P^-} \frac{1 - \cos z}{4c} q dY + \frac{1 - \cos w}{4c} p dY \right)^{\frac{1}{2}} .
\end{align*}
Thus 
\begin{equation} \label{eq:0513}
    \left| v(Q) - v (P^+)\right|    \leq h^{\frac{1}{2}} \left(\int^{P^+}_{P^-} \frac{1 - \cos z}{4c} q dY + \frac{1 - \cos w}{4c} p dY \right)^{\frac{1}{2}} .
\end{equation}
So, by (\ref{eq:0512}) and (\ref{eq:0513}) we compute that 
\begin{align*}
 |v(\tau + h ,x) - v(\tau, x)|^2 =& |v(\tau + h, x) - v(t(P^+) , x(P^+) ) + v(t(P^+) , x(P^+) ) - v(\tau, x)|^2 \\
                                  \leq& 2 \{ v(\tau + h, x) - v(t(P^+), x(P^+)) \}^2 + 2 \{ v(t(P^+), x(P^+)) - v(\tau , x)\}^2 \\
                                  \leq& 2 \{ v(Q) - v(P^+) \}^2 + 2 \{  v(P^+) - v(P) \}^2 \\
                                  \leq& 2 \left[ h^{\frac{1}{2}} \left(\int^{P^+}_{P^-} \frac{1 - \cos z}{4c} q dY + \frac{1 - \cos w}{4c} p dY \right)^{\frac{1}{2}} \right]^2\\
                                 & + 2 \left[2 \mathcal{K}^2 (\mathcal{K}h) \mu_{\tau}(]x, x+ h[) \right] \\
                                 \leq & 4 h \mu_{\tau} (]x - \mathcal{K}h, x+ \mathcal{K}h[) + 4 \mathcal{K}^3 h \mu_{\tau} (]x, x+h[)\\
                                  \leq & 4 h \mu_{\tau} (]x -\mathcal{K}h, x+ \mathcal{K}h[) (1 + \mathcal{K}^3).
\end{align*}
Thus, for all $h > 0$, 
\begin{align*}
  \| v(\tau + h, \cdot) - v(\tau, \cdot)  \|_{L^2} &= \left\{ \int |v(\tau + h, x) - v(\tau,x)|^2 dx  \right\}^{\frac{1}{2}} \\
                                                   & \leq \left\{ \int 4 (1+\mathcal{K}^3) h \mu_{\tau} (]x - \mathcal{K}h, x+ \mathcal{K}h[) \right\}^{\frac{1}{2}}\\
                                                   & \leq \left\{ 4 (\mathcal{K}^3 + 1) h^2 \mu_{\tau} (\mathbb{R}) \right\}^{\frac{1}{2}} \\
                                                   & \leq h [4 (\mathcal{K}^3 + 1) \mathcal{E}_0]^{\frac{1}{2}} \\
                                                   &\leq |\tau + h - \tau| L,
 \end{align*}
 \begin{equation} \label{eq:0514}
   \| v(\tau + h, \cdot) - v(\tau, \cdot)  \|_{L^2} \leq h [4 (\mathcal{K}^3 + 1) \mathcal{E}_0]^{\frac{1}{2}},
 \end{equation}
 where $L = [4 (\mathcal{K}^3 + 1) \mathcal{E}_0]^{\frac{1}{2}}$ is the Lipschitz constant. So, this proves the uniform Lipschitz continuous of the maps $t \mapsto v(t , \cdot)$.\qed
\section{\textbf{Regularity of trajectories}} 
In this section, we show that continuity of functions $t \mapsto v_t (t,  \cdot)$ and $t \mapsto v_x(t, \cdot)$ as functions with function value in $L^2$. It completes the proof of Theorem 1.1. \\
\indent We consider the that the initial data $ (v_0)_x$ and $v_1$ are smooth functions with compact support. In this situation, the solution $v(X,Y)$ is smooth on the $X$-$Y$ plane. Fix some time $\tau$ and denote that $\Gamma_{\tau} \coloneqq \{ (X,Y): t(X,Y) \leq \tau \}$. $\gamma_{\tau}$ is the boundary of set $\Gamma_{\tau}$. Then we claim that 
\begin{equation} \label{eq:061}
 \frac{d}{dt} v(t, \cdot) |_{t = \tau} = v_t(\tau, \cdot).
\end{equation}
By (\ref{eq:0215}), (\ref{eq:0218}), and (\ref{eq:0224}),  
\begin{equation} \label{eq:062}
\begin{split}
  v_t(\tau, x)  &\coloneqq v_X X_t + v_Y Y_t \\
  &= \frac{\sin w}{4c} p \frac{2c}{p (1+\cos w)} + \frac{\sin z}{4c} q \frac{2c}{q (1+\cos z)} \\
 &= \frac{\sin w}{2 (1+\cos w)} + \frac{\sin z}{2 (1+\cos z)}.
 \end{split}
\end{equation}     
\indent (\ref{eq:062}) define the value of $v_t(\tau, \cdot)$ at almost all the point of $x \in \mathbb{R}$. By the inequity (\ref{eq:017}) and $c(v) \geq \mathcal{K}^{-1}$, 
\begin{equation} \label{eq:065}
 \int_{\mathbb{R}} |v_t(\tau, x)|^2 dx \leq \mathcal{K}^2 \mathcal{E}(\tau) \leq \mathcal{K}^2 \mathcal{E}_0 .
\end{equation}
\indent Next, to prove (\ref{eq:061}), given $\epsilon > 0$, there exists finitely many disjoint intervals $[a_i, b_i]$ subsets of $\mathbb{R}$ with $i = 1 , 2 ... N$. We call the $A_i, B_i \in \gamma_{\tau}$ with $x(A_i) = a_i$, $x(B_i) = b_i$. Then at every point $P$ in the arcs $A_iB_i$ while $ 1 + \cos (w(P)) > \epsilon$ and  $ 1 + \cos (z(P)) > \epsilon$,
\begin{align*}
\min \{   1 + \cos (w(P)) , 1 + \cos (z(P))\} \leq 2 \epsilon .
\end{align*}
\indent We call that $J \coloneqq \bigcup_{1 \leq i \leq N} [a_i, b_i] $ as the points $P$ along the curve $\gamma_{\tau}$ that does not contain in any of the arcs $A_iB_i$ and denote that $J' \coloneqq \mathbb{R} \backslash J$. Since $v(t,x)$ is smooth in a neighbourhood of the set $\{ \tau \} \times J'$ and by the differentiability of $v$ and apply the Minkowski's inequality, 
\begin{align*}
 \lim_{h \to 0} & ~ \frac{1}{h} \left\{ \int_{\mathbb{R}} |v(\tau + h , x) - v(\tau, x) - h v_t(\tau ,x) |^p dx \right\} ^{\frac{1}{p}}\\
 \leq \lim_{ h \to 0} & ~ \frac{1}{h} \left\{ \int_j  | v(\tau + h, x) - v(\tau, x) |^p  dx \right\}^{\frac{1}{p}} + \left\{     \int_J | v_t(\tau , x)|^p dx \right\}^{\frac{1}{p}} .
\end{align*} 
\indent Now, we estimate the measure of the bad set $J$. Since $(1 +\cos w) < 2 \epsilon (1-\cos w)$ and $(1 +\cos z) < 2 \epsilon (1-\cos z)$,
\begin{align*}
meas(J) &= \int_J dx = \sum_i \int_{A_i B_i} \frac{(1+\cos w)p}{4} dX - \frac{(1+\cos z)q}{4} dY \\
        & \leq 2 \epsilon \sum_i \int_{A_i B_i} \frac{(1-\cos w)p}{4} dX - \frac{(1-\cos z)q}{4} dY \\
        & \leq 2 \epsilon \int_{\gamma_{\tau}} \frac{(1 - \cos w)p}{4} dX - \frac{(1-\cos z)q}{4} dY \\
        & \leq 2 \epsilon \mathcal{E}_0.
\end{align*}
\indent Using H\"older's inequality with exponents $\displaystyle{\frac{2}{p}}$ and $q$, we choose $\displaystyle{q = \frac{2}{2 - p}}$ so that $\displaystyle{\frac{p}{2} + \frac{1}{q} = 1}$. By (\ref{eq:0514}), 
\begin{align*}
  \int_J | v(\tau + h, x) - v(\tau, x) |^p dx & \leq meas(J)^{\frac{1}{q}} \left\{ \int_J |v(\tau, x) - v(\tau, x) |^2 dx \right\}^{\frac{p}{2}} \\
                                               & \leq [2\epsilon \mathcal{E}_0]^{\frac{1}{q}} + \left\{ h [4(\mathcal{K}^3 +1)\mathcal{E}_0]^{\frac{1}{2}} \right\}^p .                                   
\end{align*}
Thus, 
\begin{align*}
  &\lim_{h \to 0} \text{sup} \frac{1}{h} \left\{ \int_J \left|v(\tau + h, x) - v(\tau, x) \right|^p dx \right\}^{\frac{1}{p}} \\
                  \leq & \left[2\epsilon \mathcal{E}_0 \right]^{\frac{1}{pq}}  + h [4(\mathcal{K}^3 +1)\mathcal{E}_0]^{\frac{1}{2}}.
\end{align*}
Similarly, and by (\ref{eq:065}) we estimate that 
\begin{align*}
 \int_J |v_t(\tau , x) |^p dx \leq \left[meas(J)\right]^{\frac{1}{q}} \left\{ \int_J  | v_t (\tau, x) |^2dx \right\} ^{\frac{p}{2}}.
\end{align*}
Thus, 
\begin{align*}
     \left\{ \int_J |v_t(\tau , x) |^p dx  \right\}^{\frac{1}{p}} \leq [2\epsilon \mathcal{E}_0]^{\frac{1}{pq}} [\mathcal{K}^2 \mathcal{E}_0]^{\frac{1}{2}}.
\end{align*}
Since $\epsilon > 0$ is arbitrary, so we conclude that 
\begin{align*}
 \lim_{h \to 0} ~ \frac{1}{h} \left\{  \int_{\mathbb{R}}  |v(\tau + h , x) - v(\tau, x) - h v_t(\tau ,x) |^p dx \right\}^{\frac{1}{p}} = 0.
\end{align*}
\indent Next, we prove the continuity of the map $t \mapsto v_t$. 
First, we fix $\epsilon > 0$ and consider disjoint intervals $[a_i, b_i]$ subsets of $\mathbb{R}$ with $i = 1 , 2 ... N$. We call the $A_i, B_i \in \gamma_{\tau}$ with $x(A_i) = a_i$ , $x(B_i) = b_i$. Since $v$ is a smooth function on the neighbourhood of $\{ \tau\} \times J'$. By H\"older's inequality and Minkowski's inequality, we estimate that 
\begin{align*}
&\lim\sup_{ h\to 0} \int  | v_t(\tau + h, x) - v_t(\tau,x) |^p dx \\
& \leq \lim\sup_{ h \to 0}  \int_J  | v_t(\tau + h, x) - v_t(\tau,x) |^p dx     \\
& \leq \lim\sup_{ h \to 0} [meas(J)]^{\frac{1}{q}} \left\{ \int_J  |v_t( \tau + h, x) - v_t(\tau , x)|^2 dx \right\}^{\frac{q}{2}}      \\
& \leq \lim\sup_{ h \to 0} [2 \epsilon \mathcal{E}_0]^{\frac{1}{q}} \left\{ \| v_t(\tau + h, \cdot) \|_{L^2} + \| v_t(\tau, \cdot) \|_{L^2} \right\}^{\frac{q}{2}} \\
&\leq       [2 \epsilon \mathcal{E}_0]^{\frac{1}{q}}  [4 \mathcal{E}_0]^p.
\end{align*}
Since the $\epsilon > 0$ is arbitrary, so we prove the continuity. \\
 \indent For general initial data $(v_0)_x ,\, v_1 \in H^1 \cap L^4$, we consider a sequence of initial data $v^n_0 \to  v_0$, $ (v^n_0)_x \to (v_0)_x$, and $v^n_1 \to v_1$ in $\in H^1 \cap L^4$ for all $n \in \mathbb{N}$, $(v_0^n)_x$, $v_0^n$, $v_1^n \in C^{\infty}_c$. The continuity of the map $t \mapsto v_x(t, \cdot)$  with values in $L^p$ and $1 \leq p < 2$ can be proved in the same way as above.
\section{\textbf{Energy conservation}} 
\indent In this section, we provide proof of Theorem 1.3. First, we define the wave interaction potential as
\begin{equation} \label{eq:071}
  \Lambda (t) \coloneqq (\mu_{t}^- \otimes \mu_t^+) \{ (x,y): x> y \},
\end{equation}
where the $\mu_t^-$ and $\mu_t^+$ are defined in (\ref{eq:058}) and (\ref{eq:059}). And since $\mu_t^-$ and $\mu_t^+$ are absolutely continuous in Lebesgue measure, so (\ref{eq:0510}) and (\ref{eq:0511}) holds and (\ref{eq:071}) implies that 
\begin{equation} \label{eq:072}
 \Lambda (t) = \frac{1}{4} \iint_{x > y} R^2(t,x) S^2 (t,x) dx dy.
\end{equation}
\begin{lemma}
There exists a Lipschitz constant $L_0$ such that 
$$ \Lambda(t) - \Lambda(s) \leq L_0 (t - s),$$
with $t>s>0$. So the map $t \mapsto \Lambda(t)$ has bounded variation.
\end{lemma}
The Lemma is proved later in this section. \\
\indent To prove Theorem 1.3, we need to consider three sets 
\begin{align*}
& \Omega_1 \coloneqq \{  (X,Y): w(X,Y) = - \pi, \, z(X,Y) \neq - \pi, \, c'(v(X,Y)) \neq 0 \}, \\
& \Omega_2 \coloneqq \{   (X,Y): w(X,Y)  \neq \pi, \,  z(X,Y) = - \pi, \, c'(v(X,Y)) \neq 0  \} ,\\
& \Omega_3 \coloneqq \{   (X,Y): w(X,Y)= - \pi, \, z(X,Y) = -\pi,  \, c'(v(X,Y)) \neq 0 \}.
\end{align*}
From (\ref{eq:0222}), and since $w_Y \neq 0$ on $\Omega_1$ and $z_X \neq $ on $\Omega_2$, so that $meas(\Omega_1) = 0$ and $ meas(\Omega_2) = 0$. \\
We define $\Omega^*_3$ be the set of Lebesgue points of $\Omega_3$ and want to show that 
\begin{equation} \label{eq:073}
meas( \{ t(X,Y): (X,Y) \in \Omega^*_3 \} ) = 0
\end{equation}
First, we fix point $P^* \in \Omega^*_3$ and $P^* \coloneqq (X^*, Y^*)$ and claim that for $h,k > 0$,
\begin{equation} \label{eq:074}
 \lim_{h,k \to 0^+} \frac{\Lambda(\tau - h) - \Lambda(\tau + k)}{h+k} = + \infty.
\end{equation}
\indent For arbitrary $\epsilon > 0$, $\epsilon$ arbitrary small, we can find $\delta > 0$ such that for any square $Q$ with length $ l \leq \delta $ center at $P^*$, there exists a vertical segment $ \sigma $ satisfying $meas( \Omega_3 \cup \sigma  ) \geq (1 - \epsilon) l$, and a horizontal segment $\sigma'$ satisfying $meas( \Omega_3 \cup \sigma' ) \geq (1- \epsilon) l $. \\
\indent We define that 
\begin{equation} \label{eq:075}
t^+  \coloneqq \max \left\{  t(X,Y): (X,Y) \in \sigma \cup \sigma' \right\} ,
\end{equation}
\begin{equation} \label{eq:076}
t^-  \coloneqq \min \left\{  t(X,Y): (X,Y) \in \sigma \cup \sigma' \right\}.
\end{equation}
By (\ref{eq:044}), for some constant $c_0 > 0$
\begin{equation} \label{eq:077}
t^+ - t^- \leq \int_{\sigma} \frac{(1+\cos w)p}{4c} dX + \int_{\sigma'} \frac{(1+\cos z)q}{4c} dY \leq c_0 (\epsilon l)^2.
\end{equation}
(\ref{eq:077}) is Lipschitz continuous and vanished outside of a set of measure $\epsilon l$. 
Also, for some constant $c_1, c_2 > 0$,
\begin{equation} \label{eq:078}
\Lambda(t^-) - \Lambda(t^+) \geq c_1(1 - \epsilon)^2 l^2 - c_2 (t^++ - t^-).
\end{equation}
 \indent Since the choose of $ \epsilon >  0$ is arbitrary, so this implies (\ref{eq:074}). And by the Lemma 1, the map $t \mapsto \Lambda$ has bounded variation, so (\ref{eq:074}) implies (\ref{eq:073}). \\
 \indent Thus, the singular part of the $\mu_{t}$ is not trivial only if the set $\Omega_4 \coloneqq \{  P \in \gamma_{t} : w(P) = - \pi, z(P) = -\pi \}$ has positive one-dimensional measure. By the above analysis, this is restricted to a set where $c' \neq 0$ and only happens for a set of time with measure zero. \\
 {\bf{Proof of Lemma 1.}} \\
\indent From (\ref{eq:025}),
 \begin{equation} \label{eq:079}
   \left\{
       \begin{array}{ll}
            (R^2)_t - (cR^2)_x = \frac{c'}{2c} (R^2 S - S^2 R) - R(v+v^3),\\[10pt]
            (S^2)_t + (cS^2)_x = - \frac{c'}{2c} (R^2 S - S^2 R) - S(v+v^3).
        \end{array}  
        \right.
\end{equation}
 We first provide an argument valid for $v= v(t,x)$ is smooth. (\ref{eq:079}) implies that 
 \begin{align*}
 \frac{d}{dt} (4 \Lambda (t)) =& \frac{d}{dt} \iint_{x>y} R^2(t,x) S^2(t,y) dx dy \\
                               =& \iint 2 R(t,x)S^2(t,y) cR_x(t,x) + 2 S(t,y)R^2(t,x)  c S_x(t,y) \\ 
                               &+2 S(t,y)R^2(t,x) \frac{c'}{4c} (S^2(t,x) - R^2(t,x))  + 2 R(t,x)S^2(t,y)  \frac{c'}{4c} (R^2(t,x) - S^2(t,x))   \\ 
                              & -2 R(t,x)S^2(t,y) \frac{1}{2} (v(t,x) + v^3(t,x))- 2 S(t,y)R^2(t,x) \frac{1}{2} (v(t,y) + v^3(t,y)) dx dy \\                                                                                 
                               \leq& \iint c (S^2 R^2)_x + \frac{c'}{2c} (R^2 - S^2) (RS^2 - S^2 R) \\
                              &- R(t,x)S^2(t,y) [v(t,x) + v^3(t,x)]- R^2(t,x)S(t,y) [v(t,y) + v^3(t,y)] dx dy \\
                               \leq& -2 \int c R^2 S^2 dx + \int(R^2 + S^2) dx \cdot \int \frac{c'}{2c}  |R^2S -S^2 R|dx\\
                              &- \iint R(t,x)S^2(t,y) [v(t,x) + v^3(t,x)]+ R^2(t,x)S(t,y) [v(t,y) + v^3(t,y)] dx dy .
 \end{align*}
 And estimate the last term from the above calculation,
 \begin{align*}
 & \left| \iint R(t,x)S^2(t,y) [v(t,x) + v^3(t,x)] dx dy \right| \\
 & \leq \int S^2(t,y) dy \int \left|R(t,x)(v(t,x) + v^3(t,x))\right| dx \\
 &\leq \mathcal{E}_0 \| R(t,x) \|_{L^2} \| v(t,x) + v^3(t,x) \|_{L^2} \\
  & \leq \mathcal{E}_0 \mathcal{E}_0^{\frac{1}{2}} \|v(t,x) \|_{L^2} \| 1 + v(t,x) \|_{L^{\infty}}\\
  &\leq \mathcal{E}_0^2 \| 1 + v \| _{L^{\infty}}.
 \end{align*}
 Similarly,
 \begin{align*}
  & \left| \iint R^2(t,x)S(t,y) [v(t,y) + v^3(t,y)] dx dy \right| \\
  &\leq \mathcal{E}_0^2 \| 1 + v \| _{L^{\infty}}.
 \end{align*}
 Thus  
 \begin{equation}  \label{eq:0710}
  \frac{d}{dt} (4 \Lambda (t)) \leq -2 \mathcal{K}^{-1} \iint R^2 S^2 dx + 4 \mathcal{E}_0 \left\| \frac{c'}{2c} \right\|_{L^{\infty}} \int |R^2 S - S^2 R| dx + 2 \mathcal{E}_0^2 \|1 + v \|_{L^{\infty}},
 \end{equation}
where $\mathcal{K}^{-1}$ is the lower bound for the speed $c(v)$. For each $\epsilon > 0$, we have $|R| \leq \epsilon^{-\frac{1}{2}} + \epsilon^{\frac{1}{2}} R^2$. And pick any $\epsilon>0$ such that $\mathcal{K}^{-1}> 4 \mathcal{E}_0 \| \frac{c'}{2c} \|_{L^{\infty}} \cdot 2 \sqrt{\epsilon}$. \\
 Thus  
 \begin{equation}  \label{eq:0711}
   \frac{d}{dt} (4 \Lambda (t)) \leq - \mathcal{K}^{-1} \int R^2 S^2 dx + \frac{16\mathcal{E}_0^2}{\sqrt{\epsilon}} \left\|\frac{c'}{2c} \right\|_{L^{\infty}} + 2 \mathcal{E}_0^2 \|1 + v \|_{L^{\infty}}.
 \end{equation}
This yields the $L^1$ estimate:
$$ \int^{\tau}_0 \int ( | R^2 S | + | R S^2|)dx dt = \vartheta(1) \cdot [\Lambda(0) + \mathcal{E}_0^2 \tau] = \vartheta(1) \cdot (1+\tau) \mathcal{E}_0^2,$$
Where $\vartheta(1)$ is defined as a quantity and its absolute value has a uniform bound depending only on $c(v)$. Also, the map $t \mapsto \Lambda(t)$ has bounded variation on any bounded interval.
The smooth case is proved. The following provides a proof of Lemma 1 in general cases. For every $\epsilon > 0$, there exists a constant $K_{\epsilon}$ satisfying that for all $w,z$,
\begin{equation} \label{eq:0712}
\begin{split}
&|\sin z (1 - \cos w) - \sin w (1- \cos z)| \\
&\leq K_{\epsilon} \Big[\tan ^2\left(\frac{w}{2}\right)+ \tan ^2\left(\frac{z}{2}\right) \Big] (1+\cos w)(1+\cos z) + \epsilon (1-\cos w)(1-\cos z).
\end{split}
\end{equation}
For fixed $0 \leq s < t$, consider the sets $\Gamma_{s}$ and $\Gamma_{t}$ as we defined in (\ref{eq:056}) and define $\Gamma_{st} \coloneqq \Gamma_{t} \backslash \Gamma_{s}$. Recall that 
$$ dx dt = \frac{pq}{8c} (1+\cos w)(1+\cos z)dX dY.$$
We write that 
\begin{equation} \label{eq:0714}
 \int^t_s \int^{+\infty}_{-\infty} \frac{1}{4} (R^2 - S^2) dx dt  = (t-s) \mathcal{E}_0 ,
 \end{equation}
\begin{equation} \label{eq:0715} 
 \int^t_s \int^{+\infty}_{-\infty} \frac{1}{4} (R^2 - S^2) dx dt  = \iint_{\Gamma_{st}} \frac{pq}{32c} (1+\cos w)(1+\cos z)\left[\tan ^2(\frac{z}{2})+\tan ^2(\frac{w}{2})\right] dX dY.
\end{equation}
(\ref{eq:0714}) holds only on the case that $v(t,x)$ is smooth while (\ref{eq:0715}) holds for all cases.
Combine (\ref{eq:055}), (\ref{eq:058}), (\ref{eq:059}) and apply (\ref{eq:0712})-(\ref{eq:0715}), we obtain that 
\begin{align*}
\Lambda (t) - \Lambda (s) \leq& \iint_{\Gamma_{st}} \frac{(1 - \cos w)(1-\cos z)pq}{64} dX dY \\
                        &+ \mathcal{E}_0 \iint_{\Gamma_{st}} \frac{c' pq}{64c^2}[\sin z(1-\cos w) - \sin w(1-\cos z)]dX dY \\
                        &+ \mathcal{E}_0 \iint_{\Gamma_{st}} \Big\{ -\frac{pq}{32c^2} (v+v^3)(1+\cos z)\sin w - \frac{\sin w(2v^2+v^4)p}{32} \frac{c'}{8c} (\cos z - \cos w) q \\
                        &+ \frac{(1+cosw)(4v+4v^3)p}{32} \frac{1}{4c} \sin z q + \frac{(1+\cos w)}{32}\frac{c'}{8c^2}(2v^2+v^4)[\sin z-\sin w]pq \Big\} dX dY \\
                         \leq &\frac{1}{64} \iint_{\Gamma_{st}} (1-\cos w)(1-\cos z)pq dX dY \\
                        & + \mathcal{E}_0 \iint_{\Gamma_{st}} \frac{c'}{64c^2} p q \Big\{ K_{\epsilon} [\tan ^2(\frac{w}{2}) + \tan ^2(\frac{z}{2}) ] (1+\cos w)(1+\cos z)\\
                        &+ \epsilon (1-\cos w)(1-\cos z)\Big\} dX dY + \mathcal{E}_0 \iint_{\Gamma_{st}} \Big\{ -\frac{pq}{32c^2} (v+v^3)(1+\cos z)\sin w\\
                        & - \frac{\sin w(2v^2+v^4)p}{32} \frac{c'}{8c} (\cos z - \cos w) q + \frac{(1+cosw)(4v+4v^3)p}{32} \frac{1}{4c} \sin z q\\
                        & + \frac{(1+ \cos w)}{32}\frac{c'}{8c^2}(2v^2+v^4)[\sin z-\sin w]pq \Big\} dX dY \\
                       \leq & K (t - s),
\end{align*} 
for a suitable constant $K$. Thus, Lemma 1 is proved. \qed \\
\indent {\bf{Acknowledgement:}} I would like to thank Qingtian Zhang for suggesting this project and helpful discussion. I am grateful for the encouragement and support I was given throughout my thesis. Also, thanks to the Mathematics Department at UC Davis for this magnificent opportunity.


\begin{thebibliography}{99} 
\bibitem{name1} A. Bressan, and G. Chen, Generic regularity of conservative solutions to a nonlinear wave equation, \textit{Ann. I. H. Poincare-AN}, \textbf{34} (2017), 335-354.

\bibitem{name2} A. Bressan, and G. Chen, Lipschitz metrics for a class of nonlinear wave equations, \textit{Arch. Rat. Mech. Anal.}, \textbf{226} (2017), 1303-1343.

\bibitem{name3} A. Bressan, G. Chen, and Q. T. Zhang, Unique conservative solutions to a variational wave equation, \textit{Arch. Rat. Mech. Anal.}, \textbf{217} (2015), 1069-1101. 

\bibitem{name4} A. Bressan, and T. Huang, Representation of dissipative solutions to a nonlinear variational wave equation, \textit{Comm. Math. Sci.}, \textbf{14} (2016), 31-53. 

\bibitem{name5} A. Bressan, and Y. X. Zheng, Conservative solutions to a nonlinear variational wave equation, \textit{Comm. Math. Phys.} \textbf{266} (2006), no. 2, 471-497.

\bibitem{name7} G. Chen, and Y. X. Zheng, Singularity and existence for a wave system of nematic liquid crystals, \textit{J. Math. Anal. Appl.}, \textbf{398} (2013), 170-188.

\bibitem{name18} R. T. Glassey, J. K. Hunter, and Y. X. Zheng, Singularities and oscillations in a nonlinear variational wave equation, \textit{The IMA Volumes in Mathematics and its applications}, \textbf{91} (1997), 37-60.

\bibitem{newname1} H. Holden, and X. Raynaud, Global semigroup of conservative solutions of the nonlinear variational wave equation, \textit{Arch. Ration. Mech. Anal.}, \textbf{201} (2011), 871-964.

\bibitem{name8} J. C. Huang, and Y. X. Zheng, No blow-up to a variational wave equation in liquid crystals, \textit{J. Math. Phys.} \textbf{57} (2016), no. 2, 021506, 10 pp.

\bibitem{newnewname1} F. Lin, Nonlinear theory of defects in nematics liquid crystals: phase transitation and flow phenomena, \textit{Commun. Pure Appl. Math}. \textbf{42} (1989), 789-814.

\bibitem{name6} D. K. Yang, Effects of Electric Field on Liquid Crystals. \textit{Fundamentals of liquid crystal devices}. John Wiley \& Sons, 2014, 107-112. 

\bibitem{name9} P. Zhang, and Y. X. Zheng, Conservative solutions to a system of variational wave equations of nematic liquid crystals, \textit{Arch. Rat. Mech. Anal.}, \textbf{195} (2010), 701-727.

\bibitem{name10} P. Zhang, and Y. X. Zheng, Energy conservative solutions to a one-dimensional full variational wave system, \textit{Comm. Pure Appl. Math.}, \textbf{65} (2012), 683-726. 

\bibitem{name11} P. Zhang, and Y. X. Zheng, Existence and uniqueness of solutions of an asymptotic equation arising from a variational wave equation with general data, \textit{Arch. Rat. Mech. Anal.}, \textbf{155} (2000), 49-83.

\bibitem{name12} P. Zhang, and Y. X. Zheng, On the global weak solutions to a variational wave equation, \textit{Handbook of Differential Equations: Evolutionary Equations.} Vol. 2. North-Holland, 2005, 561-648.

\bibitem{name13} P. Zhang, and Y. X. Zheng, On the second-order asymptotic equation of a variational wave equation, \textit{Proc. Roy. Soc. Edinburgh}, \textbf{132} (2002), 483-509.

\bibitem{name14} P. Zhang, and Y. X. Zheng, Rarefactive solutions to a nonlinear variational wave equation of liquid crystals, \textit{Comm. Partial Differential Equations}, \textbf{26} (2001), no. 3-4, 381-419. 

\bibitem{name15} P. Zhang, and Y. X. Zheng, Singular and rarefactive solutions to a nonlinear variational wave equation, \textit{Chinese Ann. Math. Series B}, \textbf{22} (2001), 159-170.

\bibitem{name16} P. Zhang, and Y. X. Zheng, Weak solutions to a nonlinear variational wave equation, \textit{Arch. Rat. Mech. Anal.}, \textbf{166} (2003), 303-319. 
 
\bibitem{name17} P. Zhang, and Y. X. Zheng, Weak solutions to a nonlinear variational wave equation with general data, \textit{Ann. I. H. Poincare-An}, \textbf{22} (2005), 207-226.

\end{thebibliography}
\end{document}